\documentclass[preprint,12pt,nopreprintline]{elsarticle}
\usepackage{amsmath,amssymb,amsthm,mathtools,mathrsfs,bm}
\usepackage{microtype}
\usepackage{hyperref}
\biboptions{numbers,sort&compress}
\hypersetup{colorlinks=true,linkcolor=black,citecolor=black,urlcolor=blue}
\allowdisplaybreaks[3]
\numberwithin{equation}{section}
\newtheorem{theorem}{Theorem}[section]
\newtheorem{lemma}[theorem]{Lemma}

\newtheorem{corollary}[theorem]{Corollary}
\theoremstyle{remark}
\newtheorem{remark}[theorem]{Remark}
\newcommand{\R}{\mathbb R}
\newcommand{\Sn}{\mathbb S^{n-1}}
\newcommand{\Om}{\Omega}
\newcommand{\ww}{\boldsymbol w}
\newcommand{\uu}{\boldsymbol u}
\newcommand{\vv}{\boldsymbol v}
\newcommand{\zz}{\boldsymbol z}
\newcommand{\dd}{\mathrm d}
\newcommand{\prt}{\partial}
\newcommand{\norm}[1]{\left\lVert#1\right\rVert}
\newcommand{\EE}{\mathbb E}
\journal{}
\usepackage{fancyhdr}
\newcommand{\BurgersRunningTitle}{Nonradial stability of positive-far-field Burgers waves}
\fancypagestyle{burgersrunning}{%
  \fancyhf{}%
  \fancyhead[L]{\small\ifodd\value{page}\relax\itshape\BurgersRunningTitle\else\thepage\fi}%
  \fancyhead[R]{\small\ifodd\value{page}\relax\thepage\else\itshape Q. Zhao\fi}%
}
\begin{document}
\begin{frontmatter}
\title{Nonradial asymptotic stability of radial steady states for the viscous vector Burgers equation in an exterior domain with positive far field}
\author[aff1]{Qingsong Zhao}
\affiliation[aff1]{organization={School of Mathematics and Science, Nanyang Institute of Technology},
  city={Nanyang}, postcode={473004}, country={China}}
\ead{qqsszhao@nyist.edu.cn}
\begin{abstract}
We study the viscous vector Burgers equation outside a ball in $\mathbb R^n$, $n\ge2$, near a radial stationary state with positive far-field value. For every fixed $b,\mu,r_0>0$, we prove global existence and $L^\infty$ asymptotic stability for continuous, nonradial vector perturbations that vanish at infinity, satisfy the homogeneous boundary condition, and are small only in $L^\infty$. The perturbations need not be irrotational or integrable, and no initial-gradient assumption is imposed. A Bessel-weight conjugation and a divergence--curl identity yield weighted coercivity. A moving-front argument, a stopped stochastic estimate for angular momentum, and scalar Burgers comparison then recover convergence in the physical supremum norm. We also establish a quantitative Lyapunov estimate and, for perturbations compactly supported in the radial variable, an $O(t^{-1/4})$ decay bound. No additional restriction on the dimensionless parameter $br_0/\mu$ is required.
\end{abstract}
\end{frontmatter}

\section{Introduction and main result}\label{sec:intro}

\subsection{Background and relation to previous work}
Let $n\ge2$, $\mu>0$, and $r_0>0$, and denote the exterior of the ball
$B_{r_0}$ by
\[
 \Omega=\{x\in\mathbb R^n:|x|>r_0\}.
\]
We study the viscous vector Burgers system
\[
 \boldsymbol u_t+(\boldsymbol u\cdot\nabla)\boldsymbol u
 =\mu\Delta\boldsymbol u,
 \qquad (t,x)\in(0,\infty)\times\Omega,
\]
where $\boldsymbol u(t,x)\in\mathbb R^n$ and the Laplacian acts componentwise.
Writing $r=|x|$ and $\omega=x/r$, the exterior-domain problems considered
below have radial boundary and far-field data of the form
$\boldsymbol u(t,r_0\omega)=v_-\omega$ and
$\boldsymbol u(t,r\omega)\to v_+\omega$ as $r\to\infty$.
For vector fields, $|\cdot|$ is the Euclidean norm and
$\|f\|_p$ denotes the $L^p(\Omega;\mathbb R^n)$ norm; in particular,
$\|f\|_\infty=\operatorname*{ess\,sup}_{x\in\Omega}|f(x)|$.
Norms on other domains are indicated explicitly.

This system is a classical nonlinear parabolic model in which convection
and diffusion interact. Its study goes back to Burgers~\cite{Burgers1948};
for the one-dimensional scalar equation, the transformations of Hopf~\cite{Hopf1950}
and Cole~\cite{Cole1951} reduce the equation to the heat equation. Related
Cole--Hopf transformations are available for certain higher-dimensional
potential formulations~\cite{Vaganan2012}. Such reductions do not apply to a
general vector field with arbitrary angular dependence. The stability of a
radial stationary profile under nonradial vector perturbations must therefore
be studied through the full vector evolution equation.

The exterior-domain theory for multidimensional Burgers equations has been
developed mainly in the radially symmetric class. Hashimoto~\cite{Hashimoto2014,Hashimoto2016}
studied the large-time behavior of radial solutions and boundary regimes
associated with rarefaction waves. Hashimoto and Matsumura~\cite{HashimotoMatsumura2019}
analyzed convergence toward nonlinear waves, including stationary and composite
profiles. Yang, Zhao, and Zhao~\cite{YangZhaoZhao2021} established radial
asymptotic stability and decay results in higher-dimensional exterior domains;
in two dimensions, Zhao and Zhao~\cite{ZhaoZhao2021} obtained a necessary-and-sufficient
criterion for the existence of radial stationary waves and their stability under
radial perturbations. These works describe important aspects of the radial
stationary-wave branches, but do not control arbitrary angular components of a
vector perturbation.

The closest result for general multidimensional perturbations is due to
Hashimoto~\cite[Theorem~2.3]{Hashimoto2020}. That theorem treats $n\ge3$,
zero far-field state $v_+=0$, and the boundary restriction
$v_-<2\mu(n-2)/r_0$. It also assumes smallness of the stationary flow in
$L^n(\Omega)$ and of its gradient in $L^{n/2}(\Omega)$, together with
$W^{1,n/2}$-type smallness of the initial perturbation. The proof is related
to stability methods for stationary Navier--Stokes flows in exterior domains;
see Borchers and Miyakawa~\cite{BorchersMiyakawa1990} and Kozono and
Ogawa~\cite{KozonoOgawa}. Those works are methodological antecedents, rather
than direct results for the Burgers system. The present paper addresses a
different branch and a different perturbation topology: the far-field value is
strictly positive, the dimension range includes $n=2$, and smallness is required
only in the unweighted supremum norm. We do not claim an inclusion between this
$L^\infty$ data class and Hashimoto's $W^{1,n/2}$ class.

\subsection{Positive far-field stationary states}
For the positive branch, we set $v_+=b>0$, $\nu=(n-2)/2$, and $\kappa=b/(2\mu)$. Let $I_\alpha$ and $K_\alpha$ denote the modified Bessel functions of the first and second kinds, respectively. The radial stationary profile is given by
\begin{equation}\label{eq:phi}
  \boldsymbol w(x)=\phi_+(r)\omega,\qquad
  \phi_+(r)=b\frac{K_{\nu+1}(\kappa r)}{K_\nu(\kappa r)}.
\end{equation}
Its matching boundary value is $v_-:=\phi_+(r_0)$. In the notation of the
radial literature, let
\[
  V_-:=v_--\frac{\mu(n-1)}{r_0}=h(r_0),\qquad
  h(r):=\phi_+(r)-\frac{\mu(n-1)}r.
\]
For $n\ge3$, the steady-state estimates give $V_-\ge b>0$ and
$v_->2\mu(n-2)/r_0$; for $n=3$, $V_-=b$, while for $n\ge4$,
$V_->b=v_+$ and $V_-+v_+>0$. The latter sign configuration appears among
the unresolved cases in the radial classification of Yang, Zhao, and
Zhao~\cite{YangZhaoZhao2021}, although the present paper treats only the
Bessel-matched boundary value, not that entire parameter regime. Since
$w(r\omega)-b\omega\to0$ uniformly in $\omega$, the stationary state is not
in $L^n(\Omega)$; in particular, the stationary-state smallness assumption in
Hashimoto's theorem does not apply.

\subsection{Main results}
For the positive-far-field branch, fix $b>0$ and set $C_2=1$ and
$C_n=2(n-2)$ for $n\ge3$. The corresponding initial-boundary value problem is
\begin{equation}\label{eq:pde}
 \begin{cases}
 \boldsymbol u_t+(\boldsymbol u\cdot\nabla)\boldsymbol u
   =\mu\Delta\boldsymbol u,&t>0,\ x\in\Omega,\\
 \boldsymbol u(t,r_0\omega)=v_-\omega,&t\ge0,\\
 \boldsymbol u(t,r\omega)-b\omega\to0,&r\to\infty
   \text{ uniformly in }\omega,\\
 \boldsymbol u(0,x)=\boldsymbol u_0(x).&
 \end{cases}
\end{equation}
The boundary datum $v_-$ in \eqref{eq:pde} is fixed by the matching condition in \eqref{eq:phi}.

\begin{theorem}[Asymptotic stability under general nonradial perturbations]\label{thm:main}
Let $n\ge2$ and let $b,\mu,r_0>0$ be arbitrary.
Define the zero-trace space
\[
 X_0=\{f\in C(\overline\Omega;\mathbb R^n): f|_{\partial\Omega}=0,\quad
 \lim_{|x|\to\infty}f(x)=0\}.
\]
There exists $\varepsilon_*=\varepsilon_*(n,b,\mu,r_0)>0$ such that, if
\begin{equation}\label{eq:initial}
 \vv_0:=\uu_0-\ww\in X_0,\qquad \norm{\vv_0}_\infty<\varepsilon_*,
\end{equation}
then \eqref{eq:pde} admits a unique global bounded mild solution with
$\vv\in C([0,T];X_0)$ for every $T<\infty$ and
$\sup_{0<t\le\tau}\sqrt{t}\,\|\nabla\vv(s+t)\|_\infty<\infty$ on each local continuation interval $[s,s+\tau]$.
The solution is classical for positive times, and
\begin{equation}\label{eq:convergence}
 \lim_{t\to\infty}\norm{\uu(t)-\ww}_{L^\infty(\Om)}=0.
\end{equation}
No radial symmetry or irrotationality is imposed. If derivatives of the initial perturbation exist, its antisymmetric gradient need not vanish.
If, in addition, there exists $R_1>r_0$ such that $\vv_0(r,\omega)=0$ for all $r\ge R_1$
and $\omega\in\Sn$, then the following algebraic decay estimate holds:
\begin{equation}\label{eq:compactrate}
 \norm{\uu(t)-\ww}_{L^\infty(\Om)}\le C(1+t)^{-1/4},
 \qquad t\ge1.
\end{equation}
Here $C$ may depend on the support radius of the initial perturbation and on the fixed parameters.
\end{theorem}

\begin{remark}[Meaning of ``without a parameter threshold'']
This terminology means only that, for any fixed $n\ge2$ and $b,\mu,r_0>0$,
no restriction is imposed requiring $br_0/\mu$ to lie in a prescribed range.
Stability still requires the initial perturbation to have $L^\infty$ amplitude smaller than a parameter-dependent
$\varepsilon_*$. The theorem does not assert stability for large initial data and does not provide a smallness threshold uniform in
$b,\mu,r_0$; estimate~\eqref{eq:compactrate} applies only to radially compactly supported initial data and its decay exponent is not claimed to be optimal.

\end{remark}

\begin{remark}[Classical data and Lyapunov stability]
If, in addition, $\vv_0\in C_b^{2+\alpha}(\overline\Omega)$ for some $0<\alpha<1$ and
\begin{equation}\label{eq:compat}
 (\mu\Delta\vv_0-\ww\cdot\nabla\vv_0)|_{\partial\Omega}=0,
\end{equation}
standard local Dirichlet regularity also gives a classical solution up to the initial time.
These stronger conditions are not needed for the mild-solution theorem.
The following corollary makes Lyapunov stability quantitative in the supremum norm.
\end{remark}

\begin{corollary}[Quantitative Lyapunov stability]\label{cor:lyapunov}
Under the assumptions of Theorem~\ref{thm:main}, after decreasing $\varepsilon_*$ if necessary, there is a constant $C>0$ depending only on the fixed parameters such that
\begin{equation}\label{eq:lyapunov}
 \sup_{t\ge0}\|\uu(t)-\ww\|_\infty\le C\|\vv_0\|_\infty^{2/(n+2)}.
\end{equation}
\end{corollary}

\begin{remark}[Notation for vorticity in higher dimensions]
In arbitrary dimensions, ``vorticity'' is understood in its tensor form; it denotes the antisymmetric matrix
$\operatorname{Rot}(v)=(\partial_i v_j-\partial_j v_i)_{1\le i,j\le n}$.
No vanishing condition is imposed on this antisymmetric gradient anywhere in the paper.
\end{remark}

Theorem~\ref{thm:main} establishes global asymptotic stability for the positive-far-field endpoint under smallness only in the physical supremum norm. The estimate is valid for every fixed choice of $b,\mu,r_0>0$; the smallness threshold may depend on these parameters, but no additional restriction on $br_0/\mu$ is imposed. The quantitative Lyapunov estimate holds for all sufficiently small data, whereas the algebraic rate in~\eqref{eq:compactrate} is asserted only for perturbations compactly supported in the radial variable.

The proof has two analytically distinct stages. First, a Bessel-weight
conjugation of the full vector linearization leads to a divergence--curl square
identity. Together with the full Dirichlet boundary condition and exclusion of
zero modes, this gives a coercive weighted energy estimate in every dimension
$n\ge2$. The identity has the classical structure of a Witten Laplacian on
one-forms~\cite{Witten1982,HelfferNier2005}; the contribution here is its use
for the full vector Dirichlet form associated with this Burgers profile.
Second, because the conjugating Bessel weight decays at infinity, weighted
coercivity alone does not imply decay in the physical $L^\infty$ norm. We
therefore combine inner-region interpolation and a far-field supersolution
with a moving-front energy argument, a stopped stochastic estimate for angular
momentum, and signed one-dimensional Burgers comparison. Standard parabolic
boundary regularity is used to justify the positive-time estimates
~\cite{LSU,Lieberman,FurutoIwabuchi2026}; the probabilistic and scalar comparison
steps are based on standard stochastic-calculus and conservation-law tools
~\cite{RevuzYor,Serre}.

A companion manuscript by the author treats the strict negative-far-field
branch. The two studies share the Bessel stationary-state structure but have
distinct stability mechanisms; no result from the companion manuscript is used
in the present proof. With $\nu=(n-2)/2$, $\kappa=b/(2\mu)$, and
$z_0=\kappa r_0$, the present boundary value is
\[
 v_-^{(+)}=b\frac{K_{\nu+1}(z_0)}{K_\nu(z_0)}.
\]
It is exactly the critical boundary value separating the strictly negative
far-field branch, represented by a combination of $I_\nu$ and $K_\nu$, from
the positive-far-field endpoint represented by the pure $K_\nu$ branch. The
transition is singular in the global far-field topology: the two limiting
radial states have opposite far-field signs and cannot converge uniformly on
the whole exterior domain. The companion paper studies the strict negative
branch by self-adjoint spectral methods and obtains exponential decay and a
nonlinear lowest-mode expansion. The present paper instead addresses the
positive endpoint, where the decaying conjugating weight has an unbounded
reciprocal and the weighted estimate must be supplemented by transport and
angular arguments to recover physical-space asymptotic stability.

 \subsection{Notation}
We write $\mathbb S^{n-1}$ for the unit sphere and use $r=|x|$ and
$\omega=x/r$. For scalar, vector, or tensor fields $f$ on $\Omega$, we use
$|f|$ for the absolute value, Euclidean norm, or Frobenius norm,
respectively. The operator norm of a matrix is denoted by $|\cdot|_{\mathrm{op}}$.
For $1\le q<\infty$ and $q=\infty$, respectively, set
\[
 \|f\|_q=\left(\int_\Omega |f(x)|^q\,dx\right)^{1/q},
 \qquad \|f\|_\infty=\operatorname*{ess\,sup}_{x\in\Omega}|f(x)|.
\]
The domain is written
explicitly when a norm is taken on a different set. In particular,
$\|U\|_{L^p(\mathbb R)}$ refers to a scalar function on the line, and
$\|g\|_{L^2_\omega}$ is the $L^2(\mathbb S^{n-1})$ norm. We use
$\langle f,g\rangle=\int_\Omega f\cdot g\,dx$ for the real $L^2$
inner product and
\[
 \|f\|_{H^1}^2:=\|f\|_2^2+\|\nabla f\|_2^2,\qquad
 H^1=H^1(\Omega;\mathbb R^n),\quad
 H_0^1=\overline{C_c^\infty(\Omega;\mathbb R^n)}^{\,H^1}.
\]
The notation $L^p(I;X)$ denotes the usual Bochner space, and $H^{-1}$ is the dual of $H_0^1$ (with the target space understood from context). The Dirichlet Laplacian $\Delta_D$ acts componentwise. Unless explicitly fixed, $C,c>0$ denote generic constants that may change from line to line and depend on the fixed problem parameters and stationary state. Any additional dependence, for example on a support radius, a positive-time cutoff, or an initial tail, is specified where relevant. The notation $f\asymp g$ for positive functions means two-sided bounds by fixed positive multiples. The spherical differential operators are denoted by
$\nabla_S$, $\operatorname{div}_S$, and $\Delta_S$. We write $\mathbb E$ for expectation. The Bessel order and scale
are $\nu=(n-2)/2$ and $\kappa=b/(2\mu)$; $s(r)=r^{-\nu}K_\nu(\kappa r)$
is the decaying conjugating weight, whereas
$h(r)=\phi_+(r)-\mu(n-1)/r$ is the drift coefficient. We write $h_{\max}=\sup_{r\ge r_0}|h(r)|$. In the moving-front
argument of Section~\ref{sec:tails}, $c_{\mathrm{in}}$ and $\Lambda$
denote the inner and outer front speeds, respectively.

\subsection{Main ideas and organization}
The next section records the steady-state properties and constructs the global
mild solution. The subsequent sections establish the Bessel-weight energy
estimate, control the moving tails and positive-time regularity, and prove the
stochastic angular-momentum estimate and radial comparison principle. The final
section combines these ingredients to prove the main stability theorem and its
compact-support decay conclusion.
\section{Properties of the steady state and global well-posedness}\label{sec:prelim}

Set
\begin{equation}\label{eq:decomp}
 \vv=\uu-\ww=a\omega+\zz,\qquad a=\vv\cdot\omega,\quad
 \zz\cdot\omega=0,\qquad
 h(r):=\phi_+(r)-\frac{\mu(n-1)}r.
\end{equation}

\begin{lemma}[Radial coefficients of the steady state]\label{lem:steady}
For the steady state \eqref{eq:phi},
\begin{equation}\label{eq:steadyratio}
 \phi_+'+\frac{n-1}{r}\phi_+
 =\frac{\phi_+^2-b^2}{2\mu},\qquad
 0<-\phi_+'(r)\le\frac{\mu C_n}{r^2}.
\end{equation}
Moreover,
\begin{equation}\label{eq:h}
 \begin{aligned}
 \sup_{r\ge r_0}|h(r)|&=:h_{\max}<\infty,\\
 h(r)-b&=O(r^{-2}),\qquad h'(r)=O(r^{-3}).
 \end{aligned}
\end{equation}
In particular, $h\equiv b$ when $n=3$, $h\ge b$ for $n\ge3$, and $h\ge b/2$ for sufficiently large $r$ in every dimension $n\ge2$.
\end{lemma}
\begin{proof}
The radially symmetric steady-state equation is
$\phi\phi'=\mu\{\phi''+(n-1)\phi'/r-(n-1)\phi/r^2\}$.
Integrating once gives $\mu(\phi'+(n-1)\phi/r)-\phi^2/2=C_*$.
Since $\phi(r)\to b$, this identity implies that $\phi'(r)$ has the finite
limit $(C_*+b^2/2)/\mu$. A differentiable function with a finite limit
cannot have a nonzero limiting derivative, so $C_*=-b^2/2$. Hence
$\mu(\phi'+(n-1)\phi/r)=(\phi^2-b^2)/2$.
Set $\kappa=b/(2\mu)$ and $\nu=(n-2)/2$, and let
$s(r)=r^{-\nu}K_\nu(\kappa r)$. Then
$\phi_+=-2\mu(\log s)'$ yields \eqref{eq:phi} and the Riccati equation.
Conversely, the Riccati transform $\phi=-2\mu(\log s)'$ reduces the positive-far-field steady-state problem
to a modified Bessel equation. Its general solution is
\[
 s(r)=r^{-\nu}\bigl(AI_\nu(\kappa r)+BK_\nu(\kappa r)\bigr).
\]
If $A\ne0$, the exponential growth of $I_\nu(\kappa r)$ and
the exponential decay of $K_\nu(\kappa r)$ imply that the resulting branch satisfies
$\phi(r)\to-b$, which is incompatible with the positive far-field condition. Thus $A=0$,
and the positive far-field condition uniquely selects \eqref{eq:phi}; the boundary matching then determines $v_-$. 
For $\nu=(n-2)/2\ge0$, the modified Bessel function of real order has the positive integral representation
$K_\nu(x)=\int_0^\infty e^{-x\cosh s}\cosh(\nu s)\,\dd s$
(see \cite{DLMF}). By H\"older's inequality, its logarithm is convex as a function of $x>0$;
since $-\nu\log r$ is also convex and the logarithmic convexity of $K_\nu$ is strict ($(\log K_\nu)^{\prime\prime}$ is the variance of $\cosh s$ under the normalized integral weight), it follows that $\phi_+'<0$.
The integral representation above is valid for the real order $\nu\ge0$ considered here.

To obtain the dimension-dependent potential bound needed below, set
\[
 p=\frac{r\phi_+}{\mu},\qquad
 Y=-\frac{r^2\phi_+'}{\mu},\qquad t_r=\frac{br}{\mu}.
\]
Direct calculation from the Riccati equation gives
\begin{equation}\label{eq:Y}
 Y=(n-1)p-\tfrac12(p^2-t_r^2),\quad
 rp'=p-Y,\quad
 rY'=(p-n+3)Y-(n-1)p.
\end{equation}
When $n=2$, the differentiable large-argument expansion gives
\[
 \phi_+=b+\frac\mu r-\frac{\mu^2}{2br^2}+O(r^{-3}),\qquad
 Y=1-\frac{\mu}{br}+O(r^{-2}).
\]
At $Y=1$, the last equation of \eqref{eq:Y} gives $rY'=1>0$.
Since $Y<1$ near infinity, a downward crossing of $1$ would be necessary if $Y\ge1$ anywhere, which is impossible. Thus $0<Y<1$ in dimension two. The expansion yields $h-b=O(r^{-2})$, $h'=O(r^{-3})$, and eventual positivity of $h$. This proves all the assertions for $n=2$.
When $n=3$, the half-integer Bessel formula gives
$\phi_+=b+2\mu/r$, so $Y=2=C_3$ and $h=b$.
When $n\ge4$, the recurrence relation
$K_{\nu+1}(z)=K_{\nu-1}(z)+2\nu K_\nu(z)/z$
implies $p>2(n-2)$. At $Y=2(n-2)$,
the right-hand side of \eqref{eq:Y} equals $(n-3)\{p-2(n-2)\}>0$.
On the other hand, the differentiable large-argument expansion of the Bessel function \cite{DLMF} gives
\begin{equation}\label{eq:expansion}
 \begin{aligned}
 \phi_+(r)&=b+\frac{\mu(n-1)}r
 +\frac{\mu^2(n-1)(n-3)}{2br^2}+O(r^{-3}),\\
 Y(r)&\longrightarrow n-1<2(n-2).
 \end{aligned}
\end{equation}
If $Y$ reaches $2(n-2)$ at some radius, then, before tending to a lower value at infinity, it must cross this level from above to below;
however, its derivative is strictly positive at every such crossing, a contradiction.
Thus $0<Y<2(n-2)$, which gives the required bound for the potential coefficient.
Similarly, at $Y=n-1$ we have $rY'=-(n-1)(n-3)<0$;
by \eqref{eq:expansion}, $Y$ lies above $n-1$ for sufficiently large $r$, and hence it remains above this level everywhere.
It follows that $h'=\mu(n-1-Y)/r^2<0$.
Using \eqref{eq:expansion}, we obtain $h-b=O(r^{-2})$ and $h'=O(r^{-3})$;
moreover, $h$ decreases to $b$ from above. This completes the proof of \eqref{eq:h}.
\end{proof}

\begin{lemma}[Global mild solvability and preservation of the far-field condition]\label{lem:global}
For every $\vv_0\in X_0$, the problem \eqref{eq:pde} has a unique global bounded mild solution in the class of Theorem~\ref{thm:main}. It is classical for positive times and, for every $S>0$,
\begin{equation}\label{eq:C0persist}
 \vv\in C([0,S];X_0),\qquad
 \sup_{t\ge0}\|\uu(t)\|_\infty\le\sqrt n\max\{\|\uu_0\|_\infty,|v_-|,b\}.
\end{equation}
\end{lemma}
\begin{proof}
Write $C_{\rm van}(\overline\Omega)$ for continuous functions vanishing at infinity, with no boundary restriction. The Dirichlet heat evolution $S_D(t)=e^{\mu t\Delta_D}$ is strongly continuous on $X_0$ and maps $C_{\rm van}$ to $X_0$ at each positive time. It is a contraction on bounded data, and for $0<t\le1$ it satisfies
\[
 \|\nabla^j S_D(t)f\|_\infty\le C_jt^{-j/2}\|f\|_\infty,\qquad j=1,2.
\]
Indeed the killed-Brownian heat representation yields the zero boundary trace for positive times and uniform vanishing at infinity (split the data into a compactly supported part and a uniformly small tail, and bound the former by the whole-space heat kernel). The exterior-domain derivative bounds of \cite[Theorem~1.1]{FurutoIwabuchi2026} extend from smooth zero-trace data by uniform density, and to nonzero-trace inputs by splitting $S_D(t)$ into two equal time steps. Gradients vanish at infinity by interpolation using the second derivative bound. The semigroup and its gradient are continuous for positive times with values in $X_0$ and $C_{\rm van}$, respectively, by the same splitting and strong continuity on $X_0$. In particular we do not assume strong continuity at time zero on the larger space $C_{\rm van}$.

On a time interval $[s,s+\tau]$ use the Banach space of continuous $X_0$-valued functions with gradient continuous into $C_{\rm van}$ for $t>s$ and norm
\[
 \sup_{0\le t\le\tau}\|v(s+t)\|_\infty+
 \sup_{0<t\le\tau}\sqrt t\,\|\nabla v(s+t)\|_\infty.
\]
Completeness follows by identifying the uniformly convergent gradients on each positive-time subinterval with distributional derivatives of the limit. Let $G(v)=(w+v)\cdot\nabla v+(\nabla w)v$ and use the map
\[
 \Phi(v)(s+t)=S_D(t)v_s-\int_0^t S_D(t-q)G(v)(s+q)\,\dd q.
\]
Here $v_s$ is the fixed initial value at the restart time $s$. The forcing belongs to $C_{\rm van}$, but need not have zero trace. In a ball of radius $K_0$ in the above norm, the forcing is bounded by $K_0(\|w\|_\infty+K_0)q^{-1/2}+\|\nabla w\|_\infty K_0$. The difference of two such forcings is bounded by the same type of coefficient times the difference norm, using
$G(v_1)-G(v_2)=(w+v_1)\cdot\nabla(v_1-v_2)+(v_1-v_2)\cdot\nabla v_2+(\nabla w)(v_1-v_2)$.
The identities $\int_0^t q^{-1/2}\dd q=2\sqrt t$ and
$\int_0^t(t-q)^{-1/2}q^{-1/2}\dd q=\pi$ give a Lipschitz bound
$C[(\|w\|_\infty+K_0)\sqrt\tau+\|\nabla w\|_\infty\tau]$ for the integral map. Truncating the integral near each endpoint proves the required continuity by dominated convergence; the integral is $O(\sqrt t)$ in $X_0$ at the initial time. Choosing $K_0$ larger than twice the heat-evolution norm and then choosing $\tau$ small gives the contraction. The continuation time depends only on the supremum norm of the data and the fixed coefficients.

We justify positive-time classical regularity before using the maximum
principle or the weighted energy estimate. The fixed point solves the
perturbation equation in distributions, as follows by testing the heat
variation-of-constants formula against compactly supported smooth functions.
Fix $0<\delta<T$, and a bounded smooth shell $D$ whose inner boundary is the true
sphere. Let $\chi$ be a smooth spatial cutoff supported away from the outer
boundary of $D$ and equal to one on a smaller shell, and let $\eta$ be a smooth
time cutoff vanishing on $[0,\delta/4]$ and equal to one for $t\ge\delta/2$.
The field $q=\eta\chi v$ has zero initial and Dirichlet data and satisfies
\[
 q_t-\mu\Delta q=-\eta\chi G(v)+\eta'\chi v
       -2\mu\eta\nabla\chi\cdot\nabla v-\mu\eta(\Delta\chi)v.
\]
The right-hand side is bounded on the finite cylinder, since the fixed-point
gradient is bounded away from its initial time. Moreover,
\[
 q\in C([0,T];L^2(D))\cap L^2(0,T;H_0^1(D)),\qquad
 q_t\in L^2(0,T;H^{-1}(D)).
\]
The time-derivative assertion follows from the localized equation.
The linear Dirichlet heat equation with this
forcing has a strong $W_p^{2,1}$ solution for any finite $p>n+2$ by the
standard bounded-domain parabolic theory~\cite{LSU,Lieberman}. Energy uniqueness
for the zero-data difference identifies this solution with $q$. The parabolic
Sobolev embedding gives H\"older continuity of $v$ and $\nabla v$ on the smaller
shell at times $t\ge\delta/2$, with any exponent $0<\alpha<1-(n+2)/p$.
Use a second pair of spatial and time cutoffs, supported in this smaller
cylinder, with the time cutoff vanishing in a neighborhood of $\delta/2$ and
equal to one for $t\ge\delta$. The resulting heat equation has H\"older forcing
and zero initial and boundary data; its forcing vanishes near the artificial
initial face. The Dirichlet Schauder existence theorem and energy uniqueness
now give $C^{2+\alpha,1+\alpha/2}$ regularity on the remaining shell up to the true
boundary. Since $\delta$ and the shell were arbitrary, the local solution is
classical at positive times. This argument uses no global smallness assumption
and no regularity or compatibility of $v_0$ beyond $X_0$.
For each component of $u$ the equation has bounded drift $u$ on every local interval. Apply the maximum principle on truncated annuli and let the outer radius tend to infinity, using uniform far-field convergence on each compact time interval, a consequence of $C_tX_0$ continuity. This gives the componentwise bound in \eqref{eq:C0persist}. The continuation criterion therefore gives a global solution. If $v_1,v_2$ are two solutions in the asserted class, their difference $d$ obeys
\[
 d_t-\mu\Delta d+(w+v_1)\cdot\nabla d+[(\nabla w)+(\nabla v_2)]d=0.
\]
On the first common local interval Kato comparison from time $\delta>0$ gives
$\|d(t)\|_\infty\le\|d(\delta)\|_\infty
\exp(\|\nabla w\|_\infty t+\int_\delta^t\|\nabla v_2(q)\|_\infty\dd q)$.
The integral is finite uniformly as $\delta\downarrow0$, because the gradient is $O(q^{-1/2})$. Continuity at zero implies uniqueness, and restarting proves uniqueness globally. The additional classical-data assertion follows from the standard Dirichlet compatibility theorem and this uniqueness.
\end{proof}

\begin{remark}[Truncation convention for scalar comparison on the exterior domain]\label{rem:cutoffcomparison}
When comparing the Kato modulus subsolution $\rho=|\vv|$ with a positive radial supersolution
$B(t,r)$ below, we first fix an arbitrary $S<\infty$.
If $\liminf_{r\to\infty}\inf_{0\le t\le S}B(t,r)>0$,
then the uniform far-field property following \eqref{eq:C0persist} allows us to choose $R_S$ sufficiently large
so that $\rho(t,R_S,\omega)\le B(t,R_S)$ for every $t\in[0,S]$
and every $\omega\in\Sn$. We may then apply the parabolic maximum principle on a bounded annulus
and let $R_S\to\infty$.
This argument does not require a uniform far-field smallness estimate to be established in advance for all $t\ge0$.
\end{remark}

\begin{lemma}[Uniform smoothing at positive times]\label{lem:smoothing}
Let $\uu$ be the global solution from Lemma \ref{lem:global},
and set $K:=\sup_{t\ge0}\|\uu(t)\|_{L^\infty(\Om)}<\infty$.
For every $\tau_0>0$, there exist $\beta\in(0,1)$ and a constant
$C_{\tau_0}>0$, depending only on $\tau_0,n,\mu,r_0,v_-,K,\beta$, such that
\begin{equation}\label{eq:uniformsmoothing}
 \sup_{t\ge\tau_0}
 \bigl(\|\uu(t)\|_{C_b^{2+\beta}(\overline\Om)}
 +\|\vv(t)\|_{C_b^2(\overline\Om)}\bigr)
 \le C_{\tau_0}.
\end{equation}
The constant is independent of the initial derivative norms
and of the far-field truncation radius.
\end{lemma}
\begin{proof}
To prevent the large-time estimates from depending on higher-order norms of the initial data,
we first use strong $L^p$ estimates on fixed-scale positive-time parabolic cylinders and then apply Schauder estimates.
Choose $p>n+2$ and fix $0<\rho<\min\{1,r_0/8,\sqrt{\mu\tau_0}/8\}$.
All cylinders have four nested spatial scales $2\rho,\rho,\rho/2,\rho/4$,
with time scales comparable to $\rho^2$, and lie in the positive-time region $t\ge\tau_0/2$.
Writing $v=u-w$, we note that each Cartesian component of $v$ satisfies
\begin{equation}\label{eq:smoothing-linear}
 (v_i)_t-\mu\Delta v_i+u\cdot\nabla v_i
 =-\bigl((\nabla w)v\bigr)_i,
 \qquad v_i|_{\prt\Omega}=0.
\end{equation}
The drift coefficient $u$ is only required to be bounded, and the right-hand side has a uniform $L^\infty$ bound;
no a priori bound on $\nabla u$ is needed.

For interior cylinders of fixed scale, or half-cylinders obtained by flattening the true boundary
with smooth coordinate charts on the exterior sphere, local parabolic $L^p$ estimates and the absorption interpolation
inequalities for interior points and zero-Dirichlet boundary points yield
\begin{equation}\label{eq:smoothing-Lp}
 \|v\|_{W^{2,1}_p(Q_{\rho}\cap\Omega)}
 \le C_{p,\rho,K}\bigl(
 \|v\|_{L^p(Q_{2\rho}\cap\Omega)}+
 \| (\nabla w)v\|_{L^p(Q_{2\rho}\cap\Omega)}\bigr)
 \le C_{p,\rho,K}.
\end{equation}
Near the lower time face of each cylinder, we first multiply by a smooth time cutoff that vanishes there;
thus the estimate uses only local $L^p$ norms and does not involve $\|v_0\|_{C^{2+\alpha}}$.
The bounded first-order term is controlled by
$\|u\cdot\nabla v\|_{L^p}\le K\|\nabla v\|_{L^p}$
and absorbed by the local interpolation inequality.
The domain $\Omega$ is a fixed smooth exterior domain, and the boundary charts and their derivatives are uniformly bounded;
away from the boundary, we use interior cylinders of radius $\rho$ uniformly.
Consequently, $C_{p,\rho,K}$ is independent of the spatial center, the starting time of the cylinder, and the artificial truncation radius.
We use the strong-solution interior and boundary $L^p$ estimates in Chapter VII of \cite{Lieberman};
the corresponding linear parabolic $L^p$ theory in \cite{LSU} gives an alternative route.

By parabolic Sobolev embedding, after shrinking the cylinders there exists
$0<\beta<1$ (any exponent strictly less than $1-(n+2)/p$ is admissible) such that
\[
 \|v\|_{C^{\beta,\beta/2}(Q_{\rho/2}\cap\overline\Omega)}
 \le C_{\tau_0,n,\mu,r_0,K}.
\]
Since $u=w+v$ and every fixed-order derivative of $w$ is bounded, the drift $u$
and the right-hand side $-(\nabla w)v$ both have uniform parabolic H\"older norms.
Only at this stage do we apply the interior and Dirichlet-boundary Schauder estimates from Chapter IV of \cite{Lieberman}
to \eqref{eq:smoothing-linear}; after a further shrinking of the cylinders, we obtain
\[
 \|v\|_{C^{2+\beta,1+\beta/2}(Q_{\rho/4}\cap\overline\Omega)}
 \le C_{\tau_0,n,\mu,r_0,v_-,K,\beta}.
\]
The exponent $\beta$ may be decreased slightly if necessary to accommodate the boundary-coordinate changes and the embedding.
Taking the supremum over cylinders whose top time is any $t\ge\tau_0$,
and adding the fixed norms of the smooth steady field $w$, gives \eqref{eq:uniformsmoothing}.
The two-step estimate $L^p\to C^\beta\to C^{2+\beta}$ explains precisely
why the positive-time constants do not depend on higher-order norms of the initial data, and avoids
treating a merely bounded drift as a H\"older coefficient without justification.
\end{proof}

\section{Bessel-weight conjugation and global smallness}\label{sec:energy}
Set $P_R=\omega\otimes\omega$, $P_T=I-P_R$, and
\[
 s(r)=r^{-\nu}K_\nu(\kappa r),\qquad f=\log s,\quad F=\nabla f.
\]
Then $\Delta s=\kappa^2s$, $w=-2\mu F$, and $F,\nabla F\in L^\infty$.
The weight $s$ is strictly positive and bounded, with
\begin{equation}\label{eq:Wsig}
 s^2\,\dd x=W(r)\,\dd r\,\dd\omega,
 \qquad W(r)=rK_\nu(\kappa r)^2\asymp e^{-br/\mu}.
\end{equation}
In particular, $s\in L^2(\Omega)$ and $\nabla^2 f\to0$ at infinity. Thus $\|s v_0\|_2\le\|s\|_2\|v_0\|_\infty$; the weighted initial norm need not be assumed separately. This conjugation uses the same divergence--curl identity as the companion negative-far-field problem, but the weight decays rather than grows at infinity.

\begin{lemma}[Coercivity of the conjugated full vector operator]\label{lem:ground}
Let $H_+=-\mu\Delta_D+\mu \kappa^2I-2\mu\nabla^2 f$ on $L^2(\Omega;\mathbb R^n)$, with Dirichlet form $\mathcal Q_+$ and form domain $H_0^1(\Omega;\mathbb R^n)$. For every $n\ge2$ there exist $\lambda_+>0$ and $C_+<\infty$ such that
\begin{equation}\label{eq:coercivity}
 \mathcal Q_+[p]\ge\lambda_+\|p\|_2^2,
 \qquad \|p\|_{H^1}^2\le C_+\mathcal Q_+[p].
\end{equation}
\end{lemma}
\begin{proof}
The potential is bounded and symmetric, so the Dirichlet form defines a self-adjoint operator. Since $\operatorname{div}F+|F|^2=\kappa^2$, integration by parts for $p\in C_c^\infty(\Omega;\mathbb R^n)$ gives the identity
\begin{equation}\label{eq:plus-square}
 \frac{\mathcal Q_+[p]}\mu
 =\|\operatorname{div}p+F\cdot p\|_2^2
 +\sum_{i<j}\|\partial_i p_j-\partial_j p_i-F_i p_j+F_j p_i\|_2^2.
\end{equation}
Indeed the ordinary divergence and curl squares integrate to $\|\nabla p\|_2^2$; the cross terms integrate to
$\int[(\operatorname{div}F)|p|^2-2\nabla^2f[p,p]]\,\dd x$, and the squares containing $F$ alone sum to $|F|^2|p|^2$. In differential-form notation, the two terms are the squared norms of
$d_{-f}p$ and its formal adjoint $d_{-f}^*p$, where
$d_{-f}=e^f d e^{-f}=d-df\wedge$ and $p$ is identified with a one-form
\cite{Witten1982,HelfferNier2005}. We use the full zero-trace form domain,
not a relative or absolute Hodge boundary realization. Density extends the
identity to $H_0^1$. In particular the form is nonnegative.

If $\mathcal Q_+[p]=0$, the form Cauchy--Schwarz inequality implies $H_+p=0$ weakly. By the square identity, $v=s^{-1}p$ satisfies
$\partial_i v_j-\partial_j v_i=0$, $\operatorname{div}(s^2v)=0$, and $v|_{\partial\Omega}=0$ locally. Elliptic regularity for $H_+p=0$ makes these identities classical on bounded shells, up to the sphere. The full zero trace implies that
\[
 \zeta(r,\omega)=\int_{r_0}^r v(q,\omega)\cdot\omega\,\dd q
\]
has $\nabla\zeta=v$: in a local angular coordinate $\theta$, curl-freeness yields $\partial_\theta(v(q,\omega)\cdot\omega)=\partial_q(qv(q,\omega)\cdot\partial_\theta\omega)$, and the lower-end term vanishes. This construction is single-valued, including when $n=2$. Each spherical harmonic coefficient $z_\ell$ of $\zeta$ satisfies
\[
 z_\ell''+\left(\frac{n-1}{r}+2\frac{s'}s\right)z_\ell'
 -\frac{\ell(\ell+n-2)}{r^2}z_\ell=0,
 \qquad z_\ell(r_0)=z_\ell'(r_0)=0.
\]
Ordinary differential equation uniqueness gives $z_\ell=0$ for every coefficient, so $p=0$. No assumption about decay of $s^{-1}p$ at infinity is used in this step.

To prove a strictly positive gap, suppose to the contrary that $\|p_j\|_2=1$ and $\mathcal Q_+[p_j]\to0$. The bounded potential implies boundedness in $H_0^1$. After taking a subsequence, $p_j\rightharpoonup p$ in $H_0^1$ and converges strongly in $L^2$ on bounded shells. Weak lower semicontinuity of the square identity gives $\mathcal Q_+[p]=0$, hence $p=0$. Since $\nabla^2f$ is bounded and vanishes at infinity, local strong convergence and a uniform tail bound give $\int \nabla^2 f[p_j,p_j]\,\dd x\to0$. It follows that
$\mathcal Q_+[p_j]=\mu\|\nabla p_j\|_2^2+\mu \kappa^2+o(1)$, a contradiction. Thus the first inequality in \eqref{eq:coercivity} holds, and the second follows from boundedness of the potential. The argument is independent of spherical operator nonnegativity and applies in dimension two.
\end{proof}

\begin{lemma}[Nonlinear weighted energy estimate]\label{lem:weightedenergy}
There exist $M_0,\lambda>0$, depending only on the fixed parameters, such that, if $\|\vv(t)\|_\infty\le M_0$ on $[0,T]$, then
\begin{equation}\label{eq:energydecay}
 \mathcal E(t):=\frac12\int_{r_0}^\infty W(r)\|\vv(t,r,\cdot)\|_{L^2_\omega}^2\,\dd r
 \le\mathcal E(0)e^{-\lambda t},\quad 0\le t\le T.
\end{equation}
\end{lemma}
\begin{proof}
The perturbation satisfies
\begin{equation}\label{eq:pertpde}
 \vv_t+(\ww+\vv)\cdot\nabla\vv-\mu\Delta\vv=-(\nabla\ww)\vv.
\end{equation}
For $p=s\vv$, direct differentiation gives
\[
 p_t+H_+p=-(\vv\cdot\nabla)p+(F\cdot\vv)p.
\]
Notice that this expression never requires $s^{-1}$ to be bounded. At positive times $p\in H_0^1$, by uniform regularity and the exponential decay of $s$ and its derivatives. Pairing with $p$ gives
\begin{equation}\label{eq:Eidentity}
 \mathcal E'+\mathcal Q_+[p]
 =\langle -(\vv\cdot\nabla)p+(F\cdot\vv)p,p\rangle.
\end{equation}
The absolute value on the right is at most
$\|\vv\|_\infty(\|\nabla p\|_2\|p\|_2+\|F\|_\infty\|p\|_2^2)\le C_{\rm NL}\|\vv\|_\infty\mathcal Q_+[p]$
by \eqref{eq:coercivity}. Choose
\begin{equation}\label{eq:M0}
 0<M_0\le\min\{b/4,1/(4C_{\rm NL})\},\qquad \lambda=3\lambda_+/2.
\end{equation}
Then $\mathcal E'+(3/4)\mathcal Q_+[p]\le0$, giving \eqref{eq:energydecay}.
For rigor, integrate on $[\delta,T]$ with $\delta>0$, first introducing a radial cutoff equal to one for $r\le N$, zero for $r\ge2N$, and with gradient bounded by $C/N$. By Lemma~\ref{lem:smoothing}, $v$ and its derivatives are uniformly bounded there. All annular errors vanish as $N\to\infty$ by the exponential weight; at the inner sphere the trace of $p$ is zero. Finally $p\in C([0,T];L^2)$ follows from $s\in L^2$ and $v\in C([0,T];X_0)$, so letting $\delta\downarrow0$ gives the estimate from time zero. No initial gradient norm is used.
\end{proof}

\begin{lemma}[Global smallness without a parameter threshold]\label{lem:small}
There exists $\varepsilon_*>0$, depending only on $n,b,\mu,r_0$,
such that, if \eqref{eq:initial} holds, then
\begin{equation}\label{eq:uniform}
 \sup_{t\ge0}\|\vv(t)\|_\infty\le\frac{M_0}{2},\qquad
 \mathcal E(t)\le\mathcal E(0)e^{-\lambda t}.
\end{equation}
\end{lemma}
\begin{proof}
The Kato inequality for the perturbation modulus $\rho=|\vv|$ is
\begin{equation}\label{eq:kato}
 \rho_t+(\ww+\vv)\cdot\nabla\rho-\mu\Delta\rho
 -\frac{\mu C_n}{r^2}\rho\le0.
\end{equation}
Indeed, $\nabla\ww=\phi_+'P_R+\phi_+P_T/r$,
and Lemma \ref{lem:steady} gives $-\phi_+'\le\mu C_n/r^2$.
Applying the Kato regularization to $\sqrt{|\vv|^2+\delta^2}$ and letting $\delta\downarrow0$
yields \eqref{eq:kato}. In polar coordinates, after absorbing the radial
curvature term $\mu(n-1)r^{-1}\partial_r$ into
$h=\phi_+-\mu(n-1)/r$, the same inequality takes the form
\[
 \rho_t+(h+a)\rho_r+\frac1r\zz\cdot\nabla_S\rho
 -\mu\rho_{rr}-\frac{\mu}{r^2}\Delta_S\rho
 -\frac{\mu C_n}{r^2}\rho\le0.
\]
This is the radial form used when applying the scalar barriers below.

Since $h(r)\to b$, choose $R_b$ such that $h(r)\ge b/2$ for $r\ge R_b$.
First choose $A>\mu C_n/(b/2-M_0)$ and
$R_{\!*}>\max\{R_b,r_0+1,4A\}$. In this proof write $R=R_{\!*}$ and set
\begin{equation}\label{eq:barrier}
 \psi(r)=1-\frac A r,\qquad r\ge R.
\end{equation}
Then $\psi(R)>3/4$ and $\psi'>0$. If $|\vv|\le M_0$,
the inequality $h+a\ge b/2-M_0$ gives
\begin{equation}\label{eq:barrierpos}
 (h+a)\psi'-\mu\psi''-\frac{\mu C_n}{r^2}\psi
 =\frac{A(h+a)-\mu C_n}{r^2}
 +\frac{\mu A(C_n+2)}{r^3}>0.
\end{equation}

Since $\mu C_n/r^2\le\mu C_n/r_0^2$, the Kato inequality and the maximum principle yield
$\|\vv(t)\|_\infty\le e^{L_0t}\|\vv_0\|_\infty$ for $0\le t\le1$,
where $L_0:=\mu C_n/r_0^2$.

For $t\ge1$, the componentwise maximum principle first gives a uniform bound for $\uu$.
If we initially restrict $\|\vv_0\|_\infty\le1$, this bound depends only on
$n,b,\mu,r_0$. Lemma \ref{lem:smoothing}, with $\tau_0=1$, then gives

\begin{equation}\label{eq:regularity}
 \sup_{t\ge1}\bigl(\|\vv(t)\|_\infty+\|\nabla\vv(t)\|_\infty
 +\|D^2\vv(t)\|_\infty\bigr)\le K_1,
\end{equation}
where $K_1$ depends only on $n,b,\mu,r_0$ and the fixed upper bound for the initial $L^\infty$ norm,
but not on higher-order initial norms. This is precisely the independence of positive-time constants
provided by Lemma \ref{lem:smoothing}.

Suppose temporarily that $\|\vv(t)\|_\infty\le M_0$ on $[0,T]$.
By Lemma \ref{lem:weightedenergy} and $\int_{r_0}^\infty W<\infty$,
\[
 \mathcal E(t)\le C\|\vv_0\|_\infty^2e^{-\lambda t}.
\]
In the region $r\le R+1$, the weighted energy controls the usual local $L^2(\dd x)$ norm.
The uniform Lipschitz bound and a local ball/half-ball interpolation inequality in the exterior domain give
\begin{equation}\label{eq:innerinterpolation}
 \sup_{r_0\le r\le R,\omega}|\vv(t,r,\omega)|
 \le C_R\mathcal E(t)^{1/(n+2)},\qquad t\ge1.
\end{equation}
For sufficiently small initial amplitude, the energy is below $1$. We use the standard interpolation estimate
$\|f\|_{\infty,\mathcal K}\le C_{\mathcal K}
(\|f\|_{L^2(\mathcal K^+)}^{2/(n+2)}
 \|\nabla f\|_\infty^{n/(n+2)}+\|f\|_{L^2(\mathcal K^+)})$.
More precisely, on $\mathcal K=\{r_0\le r\le R\}$, choose the fixed enlarged region
$\mathcal K^+=\{r_0\le r\le R+1\}$. Since $W(r)$ has a strictly positive lower bound on this compact radial interval,
and the measures $r^{n-1}\,\dd r\,\dd\omega$
and $\dd r\,\dd\omega$ are equivalent there, we have $\|\vv(t)\|_{L^2(\mathcal K^+)}^2
\le C_R\mathcal E(t)$. Near $r=r_0$, use zero extension across the boundary of a half-ball;
the uniform bound $\|\nabla\vv(t)\|_\infty\le K_1$
then gives the stated interpolation estimate.
We may therefore choose $\|\vv_0\|_\infty<\varepsilon_*$ sufficiently small so that
\begin{equation}\label{eq:innersmall}
 \sup_{0\le t\le T,\ r_0\le r\le R}|\vv(t,r,\omega)|\le\frac{M_0}{4}.
\end{equation}

For $r\ge R$, set $B(r)=(M_0/2)\psi(r)$.
By \eqref{eq:barrierpos}, $B$ is a strict supersolution of the Kato operator.
By \eqref{eq:innersmall} and $\psi(R)>3/4$,
$|\vv(t,R,\omega)|\le M_0/4<B(R)$.
We further decrease $\varepsilon_*$ so that
$\|\vv_0\|_\infty<(M_0/2)\psi(R)$; thus the initial data are also dominated by $B$.
On every finite time interval, $\vv$ tends uniformly to zero at infinity,
whereas $B(r)\to M_0/2>0$. We may therefore apply the maximum principle
at a sufficiently distant artificial outer boundary and then let that boundary tend to infinity. It follows that
$|\vv(t,r,\omega)|\le(M_0/2)\psi(r)\le M_0/2$ for $r\ge R$.
Combining the inner and outer regions gives the strict improvement $\|\vv(t)\|_\infty\le M_0/2<M_0$.
By the continuity of $t\mapsto\vv(t)\in C_0$, a first-exit-time argument extends this estimate to every
$t\ge0$. Estimate \eqref{eq:energydecay} then also holds for all time.
\end{proof}

\begin{proof}[Proof of Corollary~\ref{cor:lyapunov}]
Write $\varepsilon=\|v_0\|_\infty\le1$; the case $\varepsilon=0$ follows by uniqueness.
Fix the barrier radius $R$ and the function $\psi$ from Lemma~\ref{lem:small},
independently of $\varepsilon$. The global smallness bound and
\eqref{eq:energydecay} give $\mathcal E(t)\le C\varepsilon^2e^{-\lambda t}$.
Combining the short-time Kato estimate with \eqref{eq:innerinterpolation} gives
\[
 \sup_{t\ge0,\,r_0\le r\le R,\,\omega\in\Sn}|v(t,r,\omega)|
 \le C_R\varepsilon^{2/(n+2)}.
\]
Choose $L\ge\max\{1,C_R\}/\psi(R)$ and put
$B_\varepsilon(r)=L\varepsilon^{2/(n+2)}\psi(r)$ for $r\ge R$.
By \eqref{eq:barrierpos} and the global bound $\|v\|_\infty\le M_0/2$,
$B_\varepsilon$ is a supersolution of the Kato inequality in the outer region.
It dominates the inner boundary data by the preceding estimate and dominates
the initial data because $\varepsilon\le\varepsilon^{2/(n+2)}$.
For any finite time interval, $B_\varepsilon$ has a strictly positive limit at
infinity and the perturbation vanishes there uniformly in time. The annular
comparison of Remark~\ref{rem:cutoffcomparison} therefore yields
$|v|\le B_\varepsilon\le L\varepsilon^{2/(n+2)}$ for every $r\ge R$.
The constants are independent of $\varepsilon$ and of the time interval.
Combining both regions proves \eqref{eq:lyapunov}.
\end{proof}

\section{Moving tails and positive-time regularity}\label{sec:tails}
By Lemma \ref{lem:small}, fix $M=M_0$ and $h_{\max}:=\sup_{r\ge r_0}|h(r)|$,
and choose
\begin{equation}\label{eq:speeds}
 0<c_{\mathrm{in}}<\min\{\lambda\mu/b,b/2\},\qquad
 \Lambda>h_{\max}+M,\qquad R\ge\max\{r_0,R_{\!*}\},
\end{equation}
where $R_{\!*}$ is the radius of the barrier constructed in Lemma~\ref{lem:small}.
Since $h_{\max}\ge b$, these choices imply $0<c_{\mathrm{in}}<\Lambda$. In particular the moving annulus has positive width $R-r_0+(\Lambda-c_{\mathrm{in}})t=O(1+t)$.

\begin{lemma}[Exponential decay behind the front]\label{lem:back}
There exist constants $C,c>0$ such that
\begin{equation}\label{eq:rear}
 \sup_{r_0\le r\le r_0+c_{\mathrm{in}} t,\omega}|\vv(t,r,\omega)|\le Ce^{-ct}.
\end{equation}
\end{lemma}
\begin{proof}
Since $W(r)\asymp e^{-br/\mu}$, the energy decays exponentially, and
$\dd x=r^{n-1}\dd r\,\dd\omega$, for $t\ge2$ and $|x|\le r_0+c_{\mathrm{in}} t$ we have
\[
 \int_{B_1(x)\cap\Om}|\vv(t,y)|^2\,\dd y
 \le C(1+t)^{n-1}e^{(b/\mu)c_{\mathrm{in}} t}\mathcal E(t)
 \le C(1+t)^{n-1}e^{-(\lambda-bc_{\mathrm{in}}/\mu)t}.
\]
By \eqref{eq:regularity} and the local Lipschitz--$L^2$ interpolation inequality,
$|\vv(t,x)|^{n+2}\le C\int_{B_1(x)\cap\Om}|\vv(t,y)|^2\,\dd y$
(provided the right-hand side is sufficiently small; the constant can be adjusted on a bounded time interval). Hence
\[
 \sup_{r\le r_0+c_{\mathrm{in}} t}|\vv(t,r,\omega)|
 \le C(1+t)^{(n-1)/(n+2)}
 e^{-(\lambda-bc_{\mathrm{in}}/\mu)t/(n+2)}.
\]
Choosing a slightly smaller exponent $c>0$ proves the claim.
\end{proof}

\begin{lemma}[Vanishing ahead of the front]\label{lem:frontC0}
For the values of $\Lambda$ and $R$ chosen in \eqref{eq:speeds},
\begin{equation}\label{eq:frontC0}
 \lim_{t\to\infty}\sup_{r\ge R+\Lambda t,\omega}|\vv(t,r,\omega)|=0.
\end{equation}
\end{lemma}
\begin{proof}
Let $\psi$ be as in \eqref{eq:barrier}, and define the radial comparison operator
$\mathscr P f=f_t+(h+a)f_r-\mu f_{rr}-\mu C_nr^{-2}f$.
By \eqref{eq:barrierpos}, $\mathscr P\psi\ge0$ for $r\ge R_{\!*}$.
Given any $\delta>0$, the condition $\vv_0\in C_0$ provides $R_\delta\ge R$ such that
$|\vv_0(r,\omega)|\le\delta\psi(r)$ for $r\ge R_\delta$.
Choose $\gamma>0$ and set
\begin{equation}\label{eq:frontbarrier}
 B_\delta(t,r)=\delta\psi(r)+C\psi(r)
 \exp\{-\gamma[r-R_\delta-(h_{\max}+M)t]+\mu\gamma^2t\}.
\end{equation}
where $C$ is chosen sufficiently large to dominate $|\vv| \le M/2$ on the initial segment $[R,R_\delta]$
and on the artificial inner boundary $r=R$.
Applying $\mathscr P$ to the second term after factoring out its positive exponential factor gives
\[
 (h+a)\psi'-\mu\psi''-\frac{\mu C_n}{r^2}\psi
 +2\mu\gamma\psi'+\gamma(h_{\max}+M-h-a)\psi\ge0.
\]
Thus $B_\delta$ is a supersolution for comparison from above.
On any finite time interval, $B_\delta(t,r)\to\delta>0$, whereas
$\vv(t,r,\omega)\to0$ uniformly in $\omega$.
We can therefore compare first at an artificial outer boundary and then let the truncation radius tend to infinity.
It follows that $|\vv|\le B_\delta$.
If $0<\gamma<(\Lambda-h_{\max}-M)/\mu$, then for $r\ge R+\Lambda t$,
\[
 |\vv(t,r,\omega)|\le\delta+C_\delta
 e^{-\gamma(\Lambda-h_{\max}-M-\mu\gamma)t}.
\]
Letting $t\to\infty$ first and then $\delta\downarrow0$ proves the claim.
\end{proof}

\begin{lemma}[Exponential decay ahead of the front for compactly supported data]\label{lem:frontCompact}
Under the assumptions of Theorem~\ref{thm:main}, suppose in addition that there exists $R_1>r_0$ such that
$\vv_0(r,\omega)=0$ for $r\ge R_1$. Then, for any fixed
$\Lambda>h_{\max}+M$ and $R\ge\max\{R_{\!*},R_1\}$, there exist $C,c>0$ such that
\begin{equation}\label{eq:frontCompact}
 \sup_{r\ge R+\Lambda t,\,\omega} |\vv(t,r,\omega)|
 \le Ce^{-ct},\qquad t\ge0.
\end{equation}
Here $M=M_0$ and $R_{\!*}$ are given by the global-smallness lemma,
and the constants $C,c$ are independent of time.
\end{lemma}
\begin{proof}
Without loss of generality, enlarge $R_1$ so that $R_1\ge R_{\!*}$.
For $r\ge R_{\!*}$, set $\psi(r)=1-A/r$.
Then $\psi\ge\psi(R_{\!*})>3/4$ and
$\mathscr P\psi\ge0$, where $\mathscr P$ is the operator in Lemma~\ref{lem:frontC0}.
Fix $\gamma>0$ and define
\[
 B(t,r)=C_0\psi(r)\exp\left\{-\gamma
 [r-R_1-(h_{\max}+M)t]+\mu\gamma^2t\right\}.
\]
The same calculation as in Lemma~\ref{lem:frontC0} gives
\[
 \frac{\mathscr PB}{C_0e^{-\gamma[r-R_1-(h_{\max}+M)t]+\mu\gamma^2t}}
 =\mathscr P\psi+2\mu\gamma\psi'
 +\gamma(h_{\max}+M-h-a)\psi\ge0.
\]
Since $\|\vv(t)\|_\infty\le M_0/2$ globally and
$\psi(R_{\!*})>3/4$, choose $C_0=M_0$. Then
on the artificial inner boundary $r=R_{\!*}$, we have $B(t,R_{\!*})\ge M_0\psi(R_{\!*})>M_0/2\ge|\vv|$ for all $t\ge0$.

At the initial time, on $[R_{\!*},R_1]$,
$B(0,r)\ge M_0\psi(r)>M_0/2\ge|\vv_0(r,\omega)|$.
whereas $\vv_0=0$ for $r\ge R_1$.
Although $B(t,r)\to0$ as $r\to\infty$,
to justify comparison on the unbounded domain we replace it, for each $\delta>0$, by
$B_\delta=B+\delta\psi$. This is still a nonnegative supersolution, and
$B_\delta(t,r)\to\delta>0$ on every finite time interval,
while $\vv(t,r,\omega)\to0$ uniformly in $\omega$.
Using the outer-boundary truncation described in Remark~\ref{rem:cutoffcomparison}
and Kato regularization, we first compare $|\vv|\le B_\delta$
and then let $\delta\downarrow0$ to obtain $|\vv|\le B$.
Finally, choose $0<\gamma<(\Lambda-h_{\max}-M)/\mu$. For
$r\ge R+\Lambda t$, we have
\[
 |\vv(t,r,\omega)|\le
 C_0e^{-\gamma(R-R_1)}
 e^{-\gamma(\Lambda-h_{\max}-M-\mu\gamma)t}.
\]
Since $R\ge R_1$, this is \eqref{eq:frontCompact}.
\end{proof}

\section{Stochastic angular momentum and a radial comparison principle}\label{sec:geometry}

\begin{lemma}[Exact equation for the radial component]\label{lem:radial}
For $a$ and $\zz$ in the decomposition \eqref{eq:decomp}, we have
\begin{equation}\label{eq:radial}
 \begin{aligned}
 a_t+(h+a)a_r+\frac1r\zz\cdot\nabla_Sa+h'a
 -\mu a_{rr}-\frac\mu{r^2}\Delta_Sa &=F_a,\\
 F_a&=-\frac{2\mu}{r^2}\operatorname{div}_S\zz+\frac{|\zz|^2}{r}.
 \end{aligned}
\end{equation}
All spherical operators here act on the angular variables on the unit sphere $\Sn$.
\end{lemma}
\begin{proof}
Since $\prt_j\omega_i=(\delta_{ij}-\omega_i\omega_j)/r$,
for a tangential vector field $\zz$ satisfying $\zz\cdot\omega=0$,
\begin{align*}
 \omega\cdot\Delta(a\omega)
 &=a_{rr}+\frac{n-1}{r}a_r
 -\frac{n-1}{r^2}a+\frac1{r^2}\Delta_Sa,\\
 \omega\cdot\Delta\zz
 &=-\frac2{r^2}\operatorname{div}_S\zz.
\end{align*}
The second identity can also be verified from
$\Delta(\omega\cdot\zz)=0$ and
$\sum_j(\prt_j\omega)\cdot\prt_j\zz=r^{-2}\operatorname{div}_S\zz$.
On the other hand, $\uu=(\phi_++a)\omega+\zz$, and
\[
 \omega\cdot[(\uu\cdot\nabla)\uu]
 =(\phi_++a)(\phi_+'+a_r)
 +r^{-1}\zz\cdot\nabla_Sa-r^{-1}|\zz|^2.
\]
Subtracting the radial steady-state equation and using
$\phi_+-\mu(n-1)/r=h$ and
$\phi_+'+\mu(n-1)/r^2=h'$, we obtain \eqref{eq:radial}.
In particular, the radial coefficient of the tangential-divergence term is $-2\mu/r^2$.
\end{proof}

Define the antisymmetric tensor
\begin{equation}\label{eq:Q}
 \begin{aligned}
 Q_{ij}(t,x)&=x_i u_j(t,x)-x_j u_i(t,x),\qquad 1\le i<j\le n,\\
 |Q|^2&=\sum_{i<j}Q_{ij}^2=r^2|\zz|^2.
 \end{aligned}
\end{equation}

\begin{lemma}[A reverse-time stochastic angular-momentum estimate]\label{lem:angular}
Let $\uu$ be a bounded solution, smooth for positive times, continuous up to time zero, and purely radial on the physical boundary,
$K:=\sup_{t\ge0,x\in\Om}|\uu(t,x)|<\infty$.
Fix $t>0$ and $x\in\Om$. Let the reverse-time characteristic be defined by
\begin{equation}\label{eq:SDE}
 \begin{aligned}
 \dd X_s&=-\uu(t-s,X_s)\,\dd s+\sqrt{2\mu}\,\dd B_s,\qquad X_0=x,\\
 \tau&=t\wedge\inf\{s\ge0:|X_s|=r_0\}.
 \end{aligned}
\end{equation}
Then
\begin{equation}\label{eq:angular}
 |Q(t,x)|\le\EE\bigl[|Q(0,X_t)|\mathbf1_{\{\tau=t\}}\bigr]
 +C_n^*K\sqrt{\mu t},
\end{equation}
where $C_n^*$ depends only on the dimension (and is distinct from the barrier constant $C_n$).
Consequently, if $\norm{Q(0)}_\infty<\infty$, then
\begin{equation}\label{eq:angular-bounded}
 \norm{Q(t)}_\infty\le\norm{Q(0)}_\infty+C_n^*K\sqrt{\mu t}.
\end{equation}
If $\vv_0\in C_0(\overline\Om)$, then for every $\sigma>0$,
\begin{equation}\label{eq:zcone}
 \sup_{r\ge\sigma t,\,\omega\in\Sn}|\zz(t,r,\omega)|\to0
 \quad(t\to\infty).
\end{equation}
\end{lemma}
\begin{proof}
On each reverse-time interval $[0,t-\delta]$, the drift is bounded and uniformly Lipschitz in space by positive-time regularity. Thus \eqref{eq:SDE} has a unique strong solution up to the first boundary hit on each such interval. Consistency and the bounded drift give a continuous extension to $s=t$ on the survival event; the time integral of the drift is absolutely convergent. The endpoint does not require differentiability of $u_0$.
Set $V_s=\uu(t-s,X_s)$. Before stopping at $\tau$, It\^o's formula gives
\begin{equation}\label{eq:martV}
 \dd V_s=\sqrt{2\mu}\,\nabla\uu(t-s,X_s)\,\dd B_s,
\end{equation}
because $-\uu_t-\uu\cdot\nabla\uu+\mu\Delta\uu=0$; see, for example, the standard It\^o and continuous-martingale theory in~\cite{RevuzYor}.
For every $0<\delta<t$, It\^o's formula up to $(t-\delta)\wedge\tau$
shows that the stopped velocity is a bounded local martingale, hence a
true martingale. These processes agree on their common time intervals.
By the continuity of $X$ up to its endpoint and the continuity of $\uu$
up to the boundary and initial time, $V_{s\wedge\tau}$ has a continuous
extension to $s=t$, held constant after a boundary hit. Bounded convergence
in the conditional-expectation identity extends the martingale property
to this endpoint. The resulting continuous bounded martingale is square
integrable, so the quadratic-variation identity gives
\begin{equation}\label{eq:bracketV}
 \sum_{j=1}^n\EE[V^j]_{\tau}
 =\EE|V_\tau|^2-|V_0|^2\le K^2.
\end{equation}
Moreover,
\[
 \EE\sup_{s\le t}|X_{s\wedge\tau}|^2
 \le C\{|x|^2+K^2t^2+\mu t\}<\infty.
\]

Strictly speaking, both $X$ and $V$ below are stopped when the process exits $\Om$;
$X_t$ appears only on the event $\{\tau=t\}$, so it is never evaluated outside the exterior domain.
For $N>\max\{r_0,|x|\}+1$ and $0<\delta<t$, define
\[
 \tau_{N,\delta}:=\tau\wedge(t-\delta)
 \wedge\inf\{s\ge0:|X_s|\ge N\}.
\]
Before this stopping time, the space-time derivatives of $\uu$ are bounded, so the following It\^o calculations
and the zero-mean property of the stochastic integrals apply directly.
Let $M^{X,i}_s=\sqrt{2\mu}B_s^i$.
For $A_s^{ij}=X_s^iV_s^j-X_s^jV_s^i$,
the finite-variation drift terms $-V^iV^j+V^jV^i$
in It\^o's product formula cancel exactly. Therefore,
\begin{equation}\label{eq:Qito}
 \begin{aligned}
 Q_{ij}(t,x)&=\EE Q_{ij}(t-\tau,X_\tau)\\
 &\quad-\EE\bigl([M^{X,i},V^j]_\tau-[M^{X,j},V^i]_\tau\bigr).
 \end{aligned}
\end{equation}
This expectation identity is first obtained at the stopping time $\tau_{N,\delta}$.
For the stopped process $V$, we have $|V|\le K$, so its components are bounded true martingales,
and $\sum_j\EE[V^j]_{\tau_{N,\delta}}\le K^2$.
Moreover, the bounded drift and the Burkholder--Davis--Gundy inequality imply
\[
 \EE\sup_{s\le t}|X_{s\wedge\tau}|^2
 \le C(n)(|x|^2+K^2t^2+\mu t).
\]
Thus the family $|X_{\tau_{N,\delta}}V_{\tau_{N,\delta}}|$
is uniformly integrable (its second moments are uniformly bounded).
On the other hand, the Kunita--Watanabe inequality gives
\begin{equation}\label{eq:KWlocal}
 \EE\bigl|[M^{X,i},V^j]_{\tau_{N,\delta}}\bigr|
 \le\sqrt{2\mu t}\,
 \bigl(\EE[V^j]_{\tau_{N,\delta}}\bigr)^{1/2}
 \le K\sqrt{2\mu t}.
\end{equation}
More explicitly, with $\rho_{N,\delta}=\tau_{N,\delta}$,
the Kunita--Watanabe inequality and Cauchy--Schwarz yield
\begin{align*}
 &\EE\left|[M^{X,i},V^j]_\tau-
 [M^{X,i},V^j]_{\rho_{N,\delta}}\right|\\
 &\quad\le
 \Bigl(2\mu\,\EE(\tau-\rho_{N,\delta})\Bigr)^{1/2}
 \Bigl(\EE\bigl([V^j]_\tau-[V^j]_{\rho_{N,\delta}}\bigr)\Bigr)^{1/2}.
\end{align*}
First let $N\to\infty$ with $\delta$ fixed, and then let $\delta\downarrow0$;
the stopping times increase to $\tau$. Moreover, $V^j$ is a bounded true martingale,
and $\EE[V^j]_\tau\le K^2$. Monotonicity of the quadratic variation and monotone convergence
show that the second factor on the right tends to zero; the first factor is always bounded by
$\sqrt{2\mu t}$ (and in fact also tends to zero).
Hence the covariations converge in $L^1$.
Furthermore, the bound $\EE\sup_{s\le t}|X_{s\wedge\tau}|^2<\infty$
together with $|V|\le K$ implies uniform integrability of the stopped products,
which converge almost surely by path continuity to $X_\tau V_\tau$.
Letting $N\to\infty$ and then $\delta\downarrow0$,
and using uniform integrability of the stopped products and the $L^1$ convergence of the covariations,
we obtain \eqref{eq:Qito}. This limiting argument does not require
any a priori higher-order smallness of $\nabla\uu$ on the full time interval.
In particular,
\begin{equation}\label{eq:KW}
 \EE\bigl|[M^{X,i},V^j]_\tau\bigr|
 \le\{\EE[M^{X,i}]_\tau\}^{1/2}
 \{\EE[V^j]_\tau\}^{1/2}
 \le K\sqrt{2\mu t}.
\end{equation}

To avoid unnecessarily enlarging the initial term when combining componentwise estimates,
take an arbitrary unit vector $A=(A_{ij})_{i<j}$ in the finite-dimensional space of antisymmetric tensors,
multiply \eqref{eq:Qito} by $A_{ij}$, and sum over $i<j$.
By Cauchy--Schwarz and \eqref{eq:KW},
the absolute value of the covariation contribution is bounded by $C_n^*K\sqrt{\mu t}$,
whereas the terminal contribution is bounded by
$\EE\{|Q(t-\tau,X_\tau)|\}$.
When $\tau<t$, we have $X_\tau=r_0\omega$ and
$V_\tau=v_-\omega$, so $Q(t-\tau,X_\tau)=0$.
When $\tau=t$, the terminal term is $Q(0,X_t)$.
Taking the supremum over all unit vectors $A$ gives \eqref{eq:angular},
and hence \eqref{eq:angular-bounded}.

Now suppose $\vv_0\in C_0$. Since the steady state is purely radial,
$|Q(0,y)|=|y\wedge\vv_0(y)|\le|y||\vv_0(y)|$. Thus, for every $\varepsilon>0$, there is a constant $C_\varepsilon$ such that
$|Q(0,y)|\le\varepsilon|y|+C_\varepsilon$, $y\in\Om$.
By \eqref{eq:SDE},
$\EE|X_\tau|\le |x|+Kt+C_n^*\sqrt{\mu t}$,
and therefore
\begin{equation}\label{eq:Qeps}
 |Q(t,x)|\le
 \varepsilon(|x|+Kt+C_n^*\sqrt{\mu t})
 +C_\varepsilon+C_n^*K\sqrt{\mu t}.
\end{equation}
Dividing by $r=|x|$ on the region $|x|\ge\sigma t$, we obtain
\[
 \limsup_{t\to\infty}\sup_{r\ge\sigma t}|\zz(t,r,\omega)|
 \le\varepsilon\bigl(1+K/\sigma\bigr).
\]
Letting $\varepsilon\downarrow0$ proves \eqref{eq:zcone}.
\end{proof}

\begin{corollary}[Global decay of the tangential component and radial forcing]\label{cor:forcing}
Fix $c_{\mathrm{in}},\Lambda,R$ as in \eqref{eq:speeds}, and set
$\mathcal A_t=\{(r,\omega):r_0+c_{\mathrm{in}} t\le r\le R+\Lambda t\}$.
Then
\begin{equation}\label{eq:force-weak}
 \norm{\zz(t)}_\infty\to0,\qquad
 \sup_{\mathcal A_t}|F_a(t,r,\omega)|=o(t^{-1}).
\end{equation}
\end{corollary}
\begin{proof}
Use \eqref{eq:zcone} for $r\ge c_{\mathrm{in}} t$ and
\eqref{eq:rear} for $r\le r_0+c_{\mathrm{in}} t$. These two regions cover the whole exterior domain,
so $\norm{\zz(t)}_\infty\to0$.

The first and second spatial derivatives of the matrix $I-\omega\otimes\omega$
are uniformly bounded for $r\ge r_0$. By \eqref{eq:regularity},
\[\sup_{t\ge1}\norm{D^2\zz(t)}_\infty\le C.\]
The uniform local interpolation inequality in the exterior domain gives
\begin{equation}\label{eq:interp}
 \norm{\nabla\zz(t)}_\infty
 \le C\Bigl(\norm{\zz(t)}_\infty^{1/2}
 \norm{D^2\zz(t)}_\infty^{1/2}+\norm{\zz(t)}_\infty\Bigr).
\end{equation}
This inequality follows from the one-dimensional Taylor formula on fixed-radius interior balls and boundary half-balls;
the covering radius depends only on $r_0$. Alternatively, it follows from a uniform extension operator on the exterior domain.
Consequently, $\norm{\nabla\zz(t)}_\infty=o(1)$.

For a tangential field $\zz\cdot\omega=0$, the spherical divergence satisfies the exact identity
$\operatorname{div}_S\zz=r\operatorname{div}_x\zz$.
By \eqref{eq:radial},
\begin{equation}\label{eq:forcebound}
 |F_a|\le \frac{C\mu}{r}|\nabla_x\zz|+\frac{|\zz|^2}{r}.
\end{equation}
On $\mathcal A_t$, we have $r\ge c_{\mathrm{in}} t$;
the preceding interpolation and decay estimates therefore imply \eqref{eq:force-weak}.
\end{proof}

\begin{corollary}[Decay of the tangential component and radial forcing for compactly supported data]\label{cor:forcingRate}
Under the assumptions of Theorem~\ref{thm:main}, suppose that $\vv_0$ is compactly supported in the radial variable.
For the moving annulus $\mathcal A_t$ defined by the fixed values of $c_{\mathrm{in}},\Lambda,R$,
there exists $C>0$ such that, for $t\ge1$,
\begin{equation}\label{eq:forceRate}
 \begin{aligned}
 \norm{\zz(t)}_\infty&\le C(1+t)^{-1/2},\\
 \sup_{\mathcal A_t}|\nabla\zz(t)|&\le C(1+t)^{-1/2},\\
 \sup_{\mathcal A_t}|F_a(t,r,\omega)|&\le Ct^{-3/2}.
 \end{aligned}
\end{equation}
\end{corollary}
\begin{proof}
Since $\ww$ is radial and $\vv_0$ is supported in $r\le R_1$,
$|Q(0,x)|=|x\wedge\vv_0(x)|\le R_1\norm{\vv_0}_\infty$.
Thus $Q(0)\in L^\infty$. By \eqref{eq:angular-bounded},
$\norm{Q(t)}_\infty\le C(1+\sqrt t)$.
On $r\ge r_0+c_{\mathrm{in}} t$, the identity $|Q|=r|\zz|$ gives
$|\zz|\le C(1+t)^{-1/2}$ for $t\ge1$.
On $r\le r_0+c_{\mathrm{in}} t$, Lemma~\ref{lem:back} gives
the same bound, and hence $\norm{\zz(t)}_\infty\le C(1+t)^{-1/2}$.
To improve the derivative estimate on the moving annulus, use $z=P_Tu$ and $\mathscr L=\partial_t+u\cdot\nabla-\mu\Delta$. Direct differentiation gives
\[
 \mathscr L z=((u\cdot\nabla)P_T-\mu\Delta P_T)u
 -2\mu\sum_{j=1}^n(\partial_j P_T)(\partial_j u).
\]
Here $|\nabla P_T|\le C/r$, $|\Delta P_T|\le C/r^2$; hence by \eqref{eq:regularity}, $|\mathscr L z|\le C/r$ at times $t\ge1$ (with the constant enlarged if necessary).
For $x$ in the moving annulus at a sufficiently large time $t$, the backward cylinder $[t-1,t]\times B_1(x)$ is interior and $|y|\ge c_{\mathrm{in}} t/2$ throughout. The global tangential bound just proved gives $\sup|z|\le Ct^{-1/2}$ on this cylinder, and $\sup|\mathscr L z|\le Ct^{-1}$. Fix $p>n+2$. The interior $W^{2,1}_p$ estimate with bounded drift, using interpolation to absorb the first-order term, on a smaller concentric cylinder followed by parabolic Sobolev embedding (as in Lemma~\ref{lem:smoothing}) gives
\[
 |\nabla z(t,x)|\le C(\sup_{[t-1,t]\times B_1(x)}|z|+\sup_{[t-1,t]\times B_1(x)}|\mathscr L z|)
 \le Ct^{-1/2}.
\]
The constant is independent of $t,x$ because the cylinder size and drift bound are fixed. For the bounded remaining time interval $t\ge1$, use \eqref{eq:regularity}. Thus, by \eqref{eq:forcebound} and $r\ge c_{\mathrm{in}} t$,
\[
 \sup_{\mathcal A_t}|F_a|\le Ct^{-1}(t^{-1/2}+t^{-1})\le Ct^{-3/2}.
\]
\end{proof}

\medskip\noindent\textbf{Preparatory estimates for signed one-dimensional comparison.}

\begin{lemma}[Signed $L^1$--$L^\infty$ and derivative estimates]\label{lem:oned}
Let $U_0\in C_c^1(\R)$, and let $U$ be the classical solution of
\begin{equation}\label{eq:oned}
 U_\tau+(b+U)U_y=\mu U_{yy},\qquad U(0,y)=U_0(y)
\end{equation}
Then, for every $\tau>0$,
\begin{align}
 \norm{U(\tau)}_{L^1(\R)}&\le\norm{U_0}_{L^1(\R)},\label{eq:L1}\\
 \norm{U(\tau)}_{L^\infty(\R)}
 &\le \sqrt{\frac{2\norm{U_0}_{L^1(\R)}}\tau},\label{eq:L1Linf}\\
 \sup_{\tau\ge0}\norm{U_y(\tau)}_\infty
 &\le \norm{U_0'}_\infty+\frac{\norm{U_0}_\infty^2}\mu.\label{eq:Ugrad}
\end{align}
If $U_0$ is nonnegative (respectively, nonpositive), then $U$ remains nonnegative (respectively, nonpositive). The constants are independent of the length of the support interval of $U_0$.
\end{lemma}
\begin{proof}
All norms of one-dimensional functions in this proof and the following
moving-annulus comparison lemma are taken on $\mathbb R$.
Define
\[
 G_0(y)=\exp\!\left[-\frac1{2\mu}\int_0^y U_0(s)\,\dd s\right]>0,
 \qquad G_\tau+bG_y=\mu G_{yy},\quad G(0)=G_0.
\]
The Cole--Hopf formula $U=-2\mu\prt_y\log G$ gives the bounded classical
solution. Since $U_0\in C_c^1(\mathbb R)$, $G_0$ is bounded above and
below by positive constants and $G_0'=-U_0G_0/(2\mu)\in L^1(\mathbb R)$.
The translated heat kernel preserves the lower bound and is contractive
on $L^1(\mathbb R)$, so
\[
 \|U(\tau)\|_{L^1(\mathbb R)}
 \le\frac{2\mu}{\inf G_0}\|G_y(\tau)\|_{L^1(\mathbb R)}
 \le\frac{2\mu}{\inf G_0}\|G_0'\|_{L^1(\mathbb R)}.
\]
In particular, $U$ is integrable uniformly on finite time intervals.
The same positive-kernel representation of $G_y/G$ gives
$\|U(\tau)\|_{L^\infty(\mathbb R)}\le\|U_0\|_{L^\infty(\mathbb R)}$.
Writing the equation in conservation form, we have
$U_\tau+\prt_y(bU+U^2/2)=\mu U_{yy}$.
The scalar Kato inequality (see \cite{Serre}) tested with
$\chi(y/R)$, where $0\le\chi\le1$ is smooth, compactly supported and
identically one near zero, yields after time integration an error bounded by
\[
 C(R^{-1}+R^{-2})\int_0^\tau\|U(s)\|_{L^1(\mathbb R)}\,\dd s.
\]
Here the flux is bounded in modulus by a constant times $|U|$.
The established integrability makes the error vanish as $R\to\infty$;
dominated convergence at both time endpoints gives \eqref{eq:L1}.
The scalar comparison principle with the constant solution zero shows that the sign is preserved if $U_0\ge0$ or $U_0\le0$.

Set $p=U_y$. Differentiating the equation gives
\begin{equation}\label{eq:p}
 p_\tau+(b+U)p_y-\mu p_{yy}=-p^2.
\end{equation}
Choose $A\ge\max\{0,\sup U_0'\}$.
The function $P(\tau)=A/(1+A\tau)$ is a spatially constant solution of \eqref{eq:p},
so $p\le A/(1+A\tau)\le1/\tau$.
If $U(\tau,y_*)=A_*>0$, then for
$y\in[y_*-A_*\tau,y_*]$ we have
$U(\tau,y)\ge A_*-(y_*-y)/\tau\ge0$,
and hence $\norm{U(\tau)}_1\ge A_*^2\tau/2$.
If $U(\tau,y_*)=-A_*<0$, then on
$[y_*,y_*+A_*\tau]$
we have $U(\tau,y)\le-A_*+(y-y_*)/\tau\le0$, which gives the same lower bound.
Approximating the supremum and infimum by sequences of points yields \eqref{eq:L1Linf}.

To obtain a uniform bound for the derivative from both sides, use the same Cole--Hopf representation.
The initial data satisfy
\[
 \left|\frac{G_0'}{G_0}\right|\le\frac{\norm{U_0}_\infty}{2\mu},\quad
 \left|\frac{G_0''}{G_0}\right|
 \le\frac{\norm{U_0'}_\infty}{2\mu}
 +\frac{\norm{U_0}_\infty^2}{4\mu^2}.
\]
Since $G$ is the convolution of the positive heat kernel with $G_0$,
the ratios $G_y/G$ and $G_{yy}/G$ are positive-weight averages of the corresponding initial ratios.
Using the identity
\[U_y=-2\mu\bigl(G_{yy}/G-(G_y/G)^2\bigr),\]
we obtain \eqref{eq:Ugrad}.
\end{proof}

\medskip\noindent\textbf{Nonradial comparison on a moving annulus.}

Fix the parameters in \eqref{eq:speeds}. For $1\le T<S$, define
\begin{equation}\label{eq:shell}
 \begin{aligned}
 \ell(t)&=r_0+c_{\mathrm{in}} t,\qquad r_+(t)=R+\Lambda t,\\
 D_{T,S}&=\{(t,r,\omega):T<t<S,\ \ell(t)<r<r_+(t)\}.
 \end{aligned}
\end{equation}

\begin{lemma}[Radial upper and lower comparison for angularly dependent solutions]\label{lem:shell}
Assume that $|a|\le M$ in $D_{T,S}$ and that
\begin{equation}\label{eq:shell-assume}
 \begin{split}
 \max\{|a(t,\ell(t),\omega)|,|a(t,r_+(t),\omega)|\}
 &\le\epsilon_{\rm side}(T),\\
 \sup_{\ell(t)\le r\le r_+(t),\omega}|F_a(t,r,\omega)|&\le f_T(t),
 \end{split}
\end{equation}
where $\epsilon_{\rm side}(T)\ge0$ and $f_T\ge0$.
There exists a constant $C$, depending on $M,c_{\mathrm{in}},\Lambda,R,r_0,n,\mu,b$ but independent of $T,S$,
and a nonnegative function $g_T\in C_c^\infty(\R)$ such that
\begin{equation}\label{eq:plateau}
 \begin{split}
 g_T(r)&\ge|a(T,r,\omega)|\quad
 (\ell(T)\le r\le r_+(T),\ \omega\in\Sn),\\
 \norm{g_T}_1&\le CT,\qquad
 \norm{g_T}_\infty+\norm{g_T'}_\infty\le C.
 \end{split}
\end{equation}
Let $U_T^\pm$ solve \eqref{eq:oned} with initial data $\pm g_T$.
Suppose that a nonnegative, nondecreasing, absolutely continuous function $E_T$ satisfies
\begin{equation}\label{eq:E}
 \begin{aligned}
 E_T(T)&\ge\epsilon_{\rm side}(T),\\
 E_T'(t)&\ge f_T(t)+Ct^{-2}+Ct^{-3}(1+E_T(t))
 \quad\text{for a.e. }t\in[T,S].
 \end{aligned}
\end{equation}
and define
$\theta_\pm(T)=0$ and $\theta_\pm'=\pm E_T$. Then
\begin{equation}\label{eq:comparison}
 U_T^-(t-T,r-\theta_-(t))-E_T(t)
 \le a(t,r,\omega)\le
 U_T^+(t-T,r-\theta_+(t))+E_T(t)
\end{equation}
for every $(t,r,\omega)\in D_{T,S}$. In particular,
\begin{equation}\label{eq:shellbound}
 \sup_{\ell(t)\le r\le r_+(t),\omega}|a(t,r,\omega)|
 \le C\sqrt{\frac T{t-T}}+\sup_{T\le s\le t}E_T(s)
 \quad(T<t\le S).
\end{equation}
\end{lemma}
\begin{proof}
On the interval $[\ell(T),r_+(T)]$, whose length is
$r_+(T)-\ell(T)=R-r_0+(\Lambda-c_{\mathrm{in}})T=O(T)$,
set $g_T=M$ and smoothly decrease it to zero on intervals of fixed length at both ends;
since $|a|\le M$, this gives \eqref{eq:plateau}.
Estimates \eqref{eq:L1Linf}--\eqref{eq:Ugrad} imply
\begin{equation}\label{eq:Ubounds}
 \begin{aligned}
 U_T^+&\ge0\ge U_T^-,\qquad
 \sup_{\tau\ge0}\norm{U_T^\pm(\tau)}_\infty\le M,\\
 \norm{\prt_yU_T^\pm}_\infty&\le C,\qquad
 \norm{U_T^\pm(t-T)}_\infty\le C\sqrt{T/(t-T)}.
 \end{aligned}
\end{equation}

Set
$\mathscr R(A):=A_t+(h+A)A_r+h'A-\mu A_{rr}$.
For the upper comparison function, which depends only on $(t,r)$,
$\overline a=U_T^+(t-T,r-\theta_+(t))+E_T(t)$,
using the equation for $U_T^+$ and $\theta_+'=E_T$, a term-by-term calculation gives
\begin{equation}\label{eq:upperres}
 \mathscr R(\overline a)-F_a
 =E_T'+(h-b)\prt_yU_T^++h'(U_T^++E_T)-F_a.
\end{equation}
Since $r\ge\ell(t)\ge c_{\mathrm{in}} t$ and $t\ge1$,
it follows from \eqref{eq:h} that
$|h-b|\le Ct^{-2}$ and $|h'|\le Ct^{-3}$.
Combining this with \eqref{eq:Ubounds} and $|F_a|\le f_T$,
we can choose the constant in \eqref{eq:E} sufficiently large, independently of $T,S$,
so that \eqref{eq:upperres} is nonnegative.
Similarly, $\underline a=U_T^-(t-T,r-\theta_-(t))-E_T(t)$
satisfies
\begin{equation}\label{eq:lowerres}
 \mathscr R(\underline a)-F_a
 =-E_T'+(h-b)\prt_yU_T^-+h'(U_T^--E_T)-F_a\le0.
\end{equation}

Although the actual solution $a$ may depend arbitrarily on the angular variables, the comparison functions $\overline a$
and $\underline a$ are independent of $\omega$.
Let $G=a-\overline a$. Subtracting the two equations gives
\begin{equation}\label{eq:G}
 \begin{aligned}
 &G_t+(h+a)G_r+\frac1r\zz\cdot\nabla_SG
 -\mu G_{rr}-\frac\mu{r^2}\Delta_SG\\
 &\qquad +(h'+\prt_r\overline a)G\le0.
 \end{aligned}
\end{equation}
By \eqref{eq:plateau}, we have $G\le0$ at $t=T$.
By the assumption in \eqref{eq:shell-assume} that \emph{each} moving lateral boundary
is separately bounded by $\epsilon_{\rm side}(T)$, together with
$U_T^+\ge0$ and $E_T\ge\epsilon_{\rm side}(T)$,
we also have $G\le0$ on both lateral boundaries.
Thus we only require the comparison functions to dominate the actual solution on the initial face and lateral boundaries of the annulus;
radial symmetry of the actual solution is not assumed.
On every finite time interval the coefficient $h'+\prt_r\overline a$ is bounded. For an absolutely continuous $E_T$, map the moving radial interval affinely onto $(0,1)$ and use the weak maximum principle: on the resulting finite cylinder the positive part of $G$ has zero trace and zero initial value, and testing the subsolution inequality with this positive part gives a Gronwall inequality for its squared $L^2$ norm. The transformed principal coefficients are uniformly elliptic on the finite window, and their spatial derivatives and the drift are bounded there. Hence $G\le0$ (for smooth $E_T$ this is also the classical maximum principle).
The same argument applied to $\underline a-a$ gives the lower comparison.
Finally, \eqref{eq:Ubounds} yields \eqref{eq:shellbound}.
\end{proof}

\begin{corollary}[A stability criterion on fixed-ratio time windows]\label{cor:window}
Assume that
\begin{equation}\label{eq:epsF}
 \epsilon_F(T):=\sup_{t\ge T}t
 \sup_{\ell(t)\le r\le r_+(t),\omega}|F_a(t,r,\omega)|\to0,
\end{equation}
and that $\epsilon_{\rm side}(T)\to0$, with $|a(t)|$ outside the moving annulus
also bounded by $\epsilon_{\rm side}(T)$ for all $t\ge T$.
Then $\norm{a(t)}_\infty\to0$.
\end{corollary}
\begin{proof}
Fix $L>1$. For $t\in[T,LT]$, choose
\[
 E_T(T)=\epsilon_{\rm side}(T),\qquad
 E_T'=\epsilon_F(T)t^{-1}+Ct^{-2}+Ct^{-3}(1+E_T).
\]
Using an integrating factor for this linear differential inequality, we obtain
\begin{equation}\label{eq:EL}
 \sup_{T\le s\le LT}E_T(s)
 \le e^{C/T^2}
 \left\{\epsilon_{\rm side}(T)+\epsilon_F(T)\log L+\frac CT\right\}.
\end{equation}
At $t=LT$, estimate \eqref{eq:shellbound}, together with the assumption outside the annulus, gives
\[
 \limsup_{T\to\infty}\norm{a(LT)}_\infty
 \le\frac C{\sqrt{L-1}}.
\]
The constant $C$ on the right is independent of $L$: the length of the plateau $g_T$ is determined only by
$T,\Lambda,c_{\mathrm{in}},R,r_0$, and the constants in \eqref{eq:L1Linf}--\eqref{eq:Ugrad}
and in the comparison remainder are independent of $S=LT$.
We first fix $L$ and let $T\to\infty$, and only then let $L\to\infty$;
no uniform-in-$L$ control of $\epsilon_F(T)\log L$ is needed.
For any sufficiently large $t$, setting $T=t/L$ in the preceding estimate controls
$\limsup_{t\to\infty}\norm{a(t)}_\infty$.
Letting $L\to\infty$ proves the claim.
\end{proof}

\section{Proof of the main theorem}\label{sec:mainproof}
\begin{proof}[Proof of Theorem \ref{thm:main}]
By Lemmas \ref{lem:global} and \ref{lem:small}, the solution exists globally,
satisfies $\|\vv(t)\|_\infty\le M_0/2$, and enjoys uniform second-order regularity at positive times.
Fix $c_{\mathrm{in}},\Lambda,R$ as in \eqref{eq:speeds}, and set
$\ell(t)=r_0+c_{\mathrm{in}} t$ and $r_+(t)=R+\Lambda t$.
Lemmas \ref{lem:back} and \ref{lem:frontC0} give
\[
 \epsilon_{\rm side}(T):=\sup_{t\ge T}\left(
 \sup_{r\le\ell(t),\omega}|\vv(t,r,\omega)|+
 \sup_{r\ge r_+(t),\omega}|\vv(t,r,\omega)|\right)\longrightarrow0.
\]
Lemma \ref{lem:angular} and Corollary \ref{cor:forcing} yield
$\|\zz(t)\|_\infty\to0$, and on the moving annulus we have
$\sup|F_a(t)|=o(t^{-1})$.
Thus both assumptions of Corollary \ref{cor:window} hold, and we obtain
$\|a(t)\|_\infty\to0$.
Since $\vv=a\omega+\zz$, this proves \eqref{eq:convergence}.

Now suppose that the initial perturbation vanishes for $r\ge R_1$. Increase $R$ so that
$R\ge\max\{R_{\!*},R_1\}$, and fix
$0<c_{\mathrm{in}}<\min\{\lambda\mu/b,b/2\}$ and $\Lambda>h_{\max}+M$.
By Lemmas~\ref{lem:back} and~\ref{lem:frontCompact},
\begin{equation}\label{eq:rateSides}
 \epsilon_{\rm side}(T)\le Ce^{-cT}.
\end{equation}
By Corollary~\ref{cor:forcingRate},
we may take $f_T(t)=Ct^{-3/2}$ on the moving annulus for $t\ge T\ge1$.

In Lemma~\ref{lem:shell}, choose
\begin{equation}\label{eq:rateE}
 E_T(T)=Ce^{-cT},\qquad
 E_T'(t)=C\bigl(t^{-3/2}+t^{-2}+t^{-3}(1+E_T(t))\bigr).
\end{equation}
Here $C$ is sufficiently large and depends only on the fixed parameters, the initial support, and the constants in the estimates;
it is independent of the right endpoint $S$ of the time window. Since
$\int_T^\infty s^{-3/2}\,\dd s=2T^{-1/2}$,
$\int_T^\infty s^{-2}\,\dd s=T^{-1}$ and
$\int_T^\infty s^{-3}\,\dd s=(2T^2)^{-1}$.
Gronwall's inequality gives, for every $S>T\ge1$,
\begin{equation}\label{eq:rateEbound}
 \sup_{T\le t\le S}E_T(t)
 \le e^{C/(2T^2)}\left(Ce^{-cT}
 +2CT^{-1/2}+CT^{-1}+CT^{-2}\right)
 \le CT^{-1/2}.
\end{equation}
Combining \eqref{eq:shellbound} with the exponential tail estimates outside the annulus, we obtain,
for $t>T\ge1$,
\begin{equation}\label{eq:rateBalance}
 \norm{\vv(t)}_\infty
 \le C\left(\sqrt{\frac{T}{t-T}}+T^{-1/2}
 +(1+t)^{-1/2}+e^{-ct}\right).
\end{equation}
Choose $T=t^{1/2}$ for sufficiently large $t$. Then
$\sqrt{T/(t-T)}=O(t^{-1/4})$ and
$T^{-1/2}=t^{-1/4}$, so
$\norm{\vv(t)}_\infty\le Ct^{-1/4}$.
Enlarging the constant on bounded time intervals gives \eqref{eq:compactrate}.
No lower bound on $br_0/\mu$ has been used.
\end{proof}
\section*{CRediT authorship contribution statement}
Qingsong Zhao: Conceptualization, Methodology, Formal analysis, Validation,
Writing -- original draft, Writing -- review \& editing.

\section*{Declaration of competing interest}
Qingsong Zhao declares that there are no known competing financial interests or personal
relationships that could have appeared to influence the work reported in this paper.

\section*{Acknowledgments}
Qingsong Zhao was supported by the National Natural Science Foundation of China under Grant No. 12401281.

\section*{Data availability}
No datasets were generated or analyzed in this theoretical study.

\section*{Declaration of generative AI and AI-assisted technologies in the manuscript preparation process}
During the preparation of this manuscript, the author used OpenAI ChatGPT to assist with language editing, literature comparison, the exploration and refinement of mathematical arguments, and manuscript revision. The author independently verified the mathematical arguments and references, reviewed and revised the AI-assisted content, and takes full responsibility for the final manuscript.


\begin{thebibliography}{99}

\bibitem{Burgers1948}
J. M. Burgers,
A mathematical model illustrating the theory of turbulence,
\emph{Adv. Appl. Mech.} \textbf{1} (1948), 171--199.
\href{https://doi.org/10.1016/S0065-2156(08)70100-5}{https://doi.org/10.1016/S0065-2156(08)70100-5}.

\bibitem{Hopf1950}
E. Hopf,
The partial differential equation $u_t+uu_x=\mu u_{xx}$,
\emph{Commun. Pure Appl. Math.} \textbf{3} (1950), 201--230.
\href{https://doi.org/10.1002/cpa.3160030302}{https://doi.org/10.1002/cpa.3160030302}.

\bibitem{Cole1951}
J. D. Cole,
On a quasi-linear parabolic equation occurring in aerodynamics,
\emph{Quart. Appl. Math.} \textbf{9} (1951), 225--236.
\href{https://doi.org/10.1090/qam/42889}{https://doi.org/10.1090/qam/42889}.

\bibitem{Vaganan2012}
B. Mayil Vaganan,
Cole--Hopf transformations for higher dimensional Burgers equations with variable coefficients,
\emph{Stud. Appl. Math.} \textbf{129} (2012), no.~3, 300--308.
\href{https://doi.org/10.1111/j.1467-9590.2012.00551.x}{https://doi.org/10.1111/j.1467-9590.2012.00551.x}.

\bibitem{Hashimoto2014}
I. Hashimoto,
Asymptotic behavior of radially symmetric solutions for the Burgers equation in several space dimensions,
\emph{Nonlinear Anal.} \textbf{100} (2014), 43--58.
\href{https://doi.org/10.1016/j.na.2014.01.004}{https://doi.org/10.1016/j.na.2014.01.004}.

\bibitem{Hashimoto2016}
I. Hashimoto,
Behavior of solutions for radially symmetric solutions for Burgers equation with a boundary corresponding to the rarefaction wave,
\emph{Osaka J. Math.} \textbf{53} (2016), 799--811.

\bibitem{HashimotoMatsumura2019}
I. Hashimoto and A. Matsumura,
Asymptotic behavior toward nonlinear waves for radially symmetric solutions of the multi-dimensional Burgers equation,
\emph{J. Differential Equations} \textbf{266} (2019), 2805--2829.
\href{https://doi.org/10.1016/j.jde.2018.08.045}{https://doi.org/10.1016/j.jde.2018.08.045}.

\bibitem{YangZhaoZhao2021}
T. Yang, H. Zhao, and Q. Zhao,
Asymptotics of radially symmetric solutions for the exterior problem of multidimensional Burgers equation (in Chinese),
\emph{Sci. Sin. Math.} \textbf{51} (2021), 1057--1072.
\href{https://doi.org/10.1360/SSM-2019-0204}{https://doi.org/10.1360/SSM-2019-0204}.

\bibitem{ZhaoZhao2021}
H. Zhao and Q. Zhao,
Radially symmetric stationary wave for two-dimensional Burgers equation,
\emph{Discrete Contin. Dyn. Syst.} \textbf{41} (2021), no.~5, 2167--2185.
\href{https://doi.org/10.3934/dcds.2020357}{https://doi.org/10.3934/dcds.2020357}.

\bibitem{Hashimoto2020}
I. Hashimoto,
Stability of the radially symmetric stationary wave of the Burgers equation with multi-dimensional initial perturbations in exterior domain,
\emph{Math. Nachr.} \textbf{293} (2020), 2348--2362.
\href{https://doi.org/10.1002/mana.201900233}{https://doi.org/10.1002/mana.201900233}.

\bibitem{BorchersMiyakawa1990}
W. Borchers and T. Miyakawa,
Algebraic $L^2$ decay for Navier--Stokes flows in exterior domains,
\emph{Acta Math.} \textbf{165} (1990), 189--227.
\href{https://doi.org/10.1007/BF02391905}{https://doi.org/10.1007/BF02391905}.

\bibitem{KozonoOgawa}
H. Kozono and T. Ogawa,
On stability of Navier--Stokes flows in exterior domains,
\emph{Arch. Rational Mech. Anal.} \textbf{128} (1994), 1--31.
\href{https://doi.org/10.1007/BF00380792}{https://doi.org/10.1007/BF00380792}.

\bibitem{Witten1982}
E. Witten,
Supersymmetry and Morse theory,
\emph{J. Differential Geom.} \textbf{17} (1982), 661--692.
\href{https://doi.org/10.4310/jdg/1214437492}{https://doi.org/10.4310/jdg/1214437492}.

\bibitem{HelfferNier2005}
B. Helffer and F. Nier,
\emph{Hypoelliptic Estimates and Spectral Theory for Fokker--Planck Operators and Witten Laplacians},
Lecture Notes in Mathematics, vol.~1862, Springer, Berlin, 2005.
\href{https://doi.org/10.1007/b104762}{https://doi.org/10.1007/b104762}.

\bibitem{LSU}
O. A. Ladyzhenskaya, V. A. Solonnikov, and N. N. Ural'tseva,
\emph{Linear and Quasilinear Equations of Parabolic Type},
Translations of Mathematical Monographs, vol. 23,
American Mathematical Society, Providence, RI, 1968.

\bibitem{Lieberman}
G. M. Lieberman,
\emph{Second Order Parabolic Differential Equations},
World Scientific, Singapore, 1996.

\bibitem{FurutoIwabuchi2026}
Y. Furuto and T. Iwabuchi,
Higher-order derivative estimates for the heat equation on a smooth domain,
\emph{Commun. Pure Appl. Anal.} \textbf{36} (2026), 55--72.
\href{https://doi.org/10.3934/cpaa.2026075}{https://doi.org/10.3934/cpaa.2026075}.

\bibitem{RevuzYor}
D. Revuz and M. Yor,
\emph{Continuous Martingales and Brownian Motion}, 3rd ed.,
Springer, Berlin, 1999.

\bibitem{Serre}
D. Serre,
\emph{Systems of Conservation Laws 1: Hyperbolicity, Entropies, Shock Waves},
Cambridge University Press, Cambridge, 1999.

\bibitem{DLMF}
NIST Digital Library of Mathematical Functions, Chapter 10: Bessel Functions,
particularly Sections 10.28, 10.29, 10.32 and 10.40,
\href{https://dlmf.nist.gov/10}{https://dlmf.nist.gov/10}.

\end{thebibliography}
\end{document}